\documentclass{amsart}

\usepackage[left=25mm,right=25mm,top=25mm,bottom=25mm]{geometry}
\usepackage{graphicx}
\usepackage[
	pdftitle={Moduli spaces of hyperbolic surfaces and their Weil-Petersson volumes},
	pdfauthor={Norman Do},
	ocgcolorlinks,
	linkcolor=linkblue,
	citecolor=linkred,
	urlcolor=linkred]
{hyperref}
\usepackage{color}
\definecolor{linkred}{rgb}{0.6,0,0}
\definecolor{linkblue}{rgb}{0,0,0.6}

\usepackage{mathpazo}
\usepackage{microtype}
\usepackage{multicol}
\usepackage{booktabs}

\theoremstyle{plain}
\newtheorem{theorem}{Theorem}

\newtheorem{proposition}{Proposition}[section]

\newtheorem{corollary}[proposition]{Corollary}

\theoremstyle{definition}
\newtheorem{example}[proposition]{Example}
\newtheorem{definition}[proposition]{Definition}
\newtheorem{remark}[proposition]{Remark}

\newcommand{\M}{\overline{\mathcal M}}

\newcommand{\bb}{{\mathbf b}}
\newcommand{\N}{\overline{N}}
\newcommand{\PP}{\mathbb{P}}

\newcommand{\bc}{\mathbb{C}}
\newcommand{\bn}{\mathbb{N}}
\newcommand{\bp}{\mathbb{P}}
\newcommand{\bq}{\mathbb{Q}}
\newcommand{\bz}{\mathbb{Z}}
\newcommand{\f}{\mathcal{F}}
\newcommand{\ci}{\mathcal{I}}
\newcommand{\y}{\mathcal{Y}}
\newcommand{\z}{\mathcal{Z}}
\newcommand{\modm}{\mathcal{M}}
\newcommand{\fat}{\f\hspace{-.3mm}{\rm at}_{g,n}}
\newcommand{\fats}{\f\hspace{-.3mm}{\rm at}^{\rm stable}_{g,n}}
\newcommand{\fato}{\f\hspace{-.3mm}{\rm at}}
\newcommand{\ie}{i.e.\ }

\begin{document}

\title{Counting lattice points in compactified moduli spaces of curves}
\author{Norman Do \and Paul Norbury}
\address{Department of Mathematics and Statistics, The University of Melbourne, Victoria 3010, Australia}
\email{\href{mailto:normdo@gmail.com}{normdo@gmail.com}, \href{mailto:pnorbury@ms.unimelb.edu.au}{pnorbury@ms.unimelb.edu.au}}

\date{\today}
\begin{abstract}
We define and count lattice points in the moduli space $\M_{g,n}$ of stable genus $g$ curves with $n$ labeled points. This extends a construction of the second author for the uncompactified moduli space ${\mathcal M}_{g,n}$. The enumeration produces polynomials whose top degree coefficients are tautological intersection numbers on $\M_{g,n}$ and whose constant term is the orbifold Euler characteristic of $\M_{g,n}$. We prove a recursive formula which can be used to effectively calculate these polynomials. One consequence of these results is a simple recursion relation for the orbifold Euler characteristic of $\M_{g,n}$.
\end{abstract}

\maketitle

\tableofcontents

\section{Introduction} \label{introduction}

Lattice points in the moduli space ${\mathcal M}_{g,n}$ of smooth genus $g$ curves with $n$ labeled points were defined and counted in \cite{NorCou}. For positive integers $b_1, b_2, \ldots, b_n$, define $\z_{g,n}(b_1, b_2, \ldots,b_n) \subset \modm_{g,n}$ to consist of any smooth curve $\Sigma$ with labeled points $(p_1, p_2, \ldots, p_n)$ that possesses a morphism $f: \Sigma \to \bp^1$ satisfying the following three conditions.
\begin{enumerate}
\item[(C1)] \label{p1} $f$ has degree $b_1 + b_2 + \cdots + b_n$ and is regular over $\bp^1 \setminus \{0,1,\infty\}$.
\item[(C2)] \label{p2} $f^{-1}(\infty)= \{p_1, p_2,..., p_n\}$ with ramification $b_k$ at $p_k$. Each point in $f^{-1}(1)$ has ramification of order 2.
\item[(C3)] \label{p3} There are no points with ramification of order 1 over $0 \in \bp^1$.
\end{enumerate}
We count the number of points in the finite set $\z_{g,n}(b_1, b_2, \ldots, b_n)$ taking into account the orbifold nature of $\modm_{g,n}$. More precisely, a point $\Sigma \in \z_{g,n}(b_1, b_2, \ldots, b_n)$ is counted with weight equal to the reciprocal of the order of its automorphism group. The weighted count is conveniently expressed by the orbifold Euler characteristic of $\z_{g,n}(b_1, b_2, \ldots, b_n)$.

\begin{definition} \label{th:defngn}
For positive integers $b_1, b_2, \ldots, b_n$, define
\[
N_{g,n}(b_1, b_2, \ldots, b_n)= \chi\left(\z_{g,n}(b_1, b_2, \ldots, b_n)\right)\in\bq.
\]
\end{definition}
It was shown in \cite{NorCou} that $N_{g,n}(b_1, b_2, \ldots,b_n)$ is recursively calculable and quasi-polynomial in $b_1^2, b_2^2, \ldots, b_n^2$ in the sense that it is polynomial on each coset of the sublattice $2\bz^n \subset \bz^n$.

In this paper, we propose a lattice point count $\N_{g,n}(b_1, b_2, \ldots, b_n)$ which augments $N_{g,n}(b_1, b_2, \ldots, b_n)$ in a natural way. The extra contribution arises from stable genus $g$ curves with $n$ labeled points in the boundary divisor of the Deligne--Mumford compactification $\M_{g,n}$. Recall that an algebraic curve is called stable if its singularities are nodal and its automorphism group is finite. As above, for positive integers $b_1, b_2, \ldots, b_n$, define $\overline{\z}_{g,n}(b_1, b_2, \ldots, b_n) \subset \M_{g,n}$ to consist of any stable curve $\Sigma$ with labeled points $(p_1, p_2, \ldots, p_n)$ that possesses a morphism $f:\Sigma\to\bp^1$ satisfying conditions~(C1) and (C2) above as well as the following.
\begin{enumerate}
\item[(C3')] Every point with ramification of order 1 over $0 \in \bp^1$ is a node.
\end{enumerate}
Nodes and ghost components --- irreducible components without labeled points --- necessarily lie in the fibre over $0 \in \bp^1$. The set $\overline{\z}_{g,n}(b_1, b_2, \ldots, b_n)$ is no longer finite since ghost components can introduce moduli. Nevertheless, we can generalise the definition above and virtually count points in $\M_{g,n}$ using the orbifold Euler characteristic.
\begin{definition}
For positive integers $b_1, b_2, \ldots, b_n$, define
\[
\overline{N}_{g,n}(b_1, b_2, \ldots, b_n) = \chi\left(\overline{\z}_{g,n}(b_1, b_2, \ldots, b_n)\right)\in\bq.
\]
\end{definition}

\begin{remark} \label{th:suborb}
Given $(b_1, b_2, \ldots, b_n)\in\bz_+^n$ if a stable curve admits a morphism satisfying (C1), (C2) and (C3') then that morphism is unique, and hence it makes sense to write $\overline{\z}_{g,n}(b_1, b_2, \ldots, b_n)$ as a subset of $\overline{\modm}_{g,n}$. Furthermore, any automorphism of a curve in $\overline{\z}_{g,n}(b_1, b_2, \ldots, b_n)$ fixes its morphism satisfying (C1), (C2) and (C3'), \ie the two automorphism groups coincide. So $\overline{\z}_{g,n}(b_1, b_2, \ldots, b_n)$ is naturally a suborbifold of $\overline{\modm}_{g,n}$. See Section~\ref{sec:fatgraph} for more details.
\end{remark}

The compactified lattice point count $\N_{g,n}(b_1, b_2, \ldots, b_n)$ has a particularly nice structure, as evidenced by the following result which is an analogue of results concerning the uncompactified count $N_{g,n}(b_1, b_2, \ldots, b_n)$, \cite{NorCou}.
\begin{theorem} \label{th:nbarprop} ~
\begin{itemize}
\item The compactified lattice point count $\N_{g,n}(b_1, b_2, \ldots, b_n)$ is a symmetric quasi-polynomial in $b_1^2, b_2^2, \ldots, b_n^2$ of degree $3g-3+n$ in the sense that it is polynomial on each coset of the sublattice $2\bz^n \subset \bz^n$.

\item If $\alpha_1 + \alpha_2 + \cdots + \alpha_n = 3g-3+n$, then the coefficient of $b_1^{2\alpha_1} b_2^{2\alpha_2} \cdots b_n^{2\alpha_n}$ in $\N_{g,n}(b_1, b_2, \ldots, b_n)$ is the following intersection number of psi-classes $\psi_1, \psi_2, \ldots, \psi_n \in H^2(\M_{g,n}; \bq)$.
\[
\frac{1}{2^{5g-6+2n} \alpha_1! \alpha_2! \cdots \alpha_n!} \int_{\M_{g,n}} \psi_1^{\alpha_1} \psi_2^{\alpha_2} \cdots \psi_n^{\alpha_n}
\]

\item The constant coefficient of $\N_{g,n}(b_1, b_2, \ldots, b_n)$ is the orbifold Euler characteristic of $\M_{g,n}$.
\[
\N_{g,n}(0, 0, \ldots, 0) = \chi(\M_{g,n})
\]
\end{itemize}
\end{theorem}
The polynomials on each coset of the sublattice $2\bz^n \subset \bz^n$ that represent $\N_{g,n}(b_1, b_2, \ldots, b_n)$ are denoted $\N_{g,n}^{(k)}(b_1, b_2, \ldots, b_n)$ where $k$ is the number of odd $b_i$. Note that the enumeration $\N_{g,n}(b_1, b_2, \ldots, b_n)$ is defined only when $b_1, b_2, \ldots, b_n$ are positive integers. However, its quasi-polynomial behaviour allows us to evaluate $N_{g,n}(b_1, b_2, \ldots, b_n)$ for arbitrary integers $b_1, b_2, \ldots, b_n$.

The tautological intersection numbers stored in the top degree coefficients of $\N_{g,n}(\bb)$ are precisely those which are governed by the Witten--Kontsevich theorem \cite{KonInt, WitTwo}. The orbifold Euler characteristic for the Deligne--Mumford compactification $\M_{g,n}$ is computed in \cite{BHaEul}, though not in explicit form. It is interesting that these two calculations should appear together in the context of counting lattice points in $\M_{g,n}$. We remark that it is currently unknown whether or not the intermediate coefficients of $\N_{g,n}(\bb)$ store topological information about $\M_{g,n}$.

The following recursive formula can be used to effectively compute $\N_{g,n}(\bb)$ from the base cases $\N_{0,3}(b_1, b_2, b_3)$ and $\N_{1,1}(b_1)$.

\begin{theorem} \label{th:recursion}
Let $S = \{1,2, 3, \ldots, n\}$ and for an index set $I = \{i_1, i_2, \ldots, i_m\}$, let $\bb_I = (b_{i_1}, b_{i_2}, \ldots, b_{i_m})$. The compactified lattice point count satisfies the following recursive formula, 
\begin{align} \label{eq:rec}
\left(\sum_{i=1}^nb_i\right) \N_{g,n}(\bb_S) &= \sum_{i\neq j} \sum_{p+q=b_i+b_j} f(p)q \N_{g,n-1}(p, \bb_{S \setminus \{i,j\}})\\
+ \frac{1}{2}\sum_i&\sum_{p+q+r=b_i} f(p) f(q) r \biggl[\N_{g-1,n+1}(p, q, \bb_{S\setminus \{i\}})
+ \hspace{-5mm}\sum_{\substack{g_1+g_2=g\\I_1 \sqcup I_2 = S\setminus \{i\}}} \hspace{-3mm}\N_{g_1,|I_1|+1}(p, \bb_{I_1}) \N_{g_2,|I_2|+1}(q, \bb_{I_2}) \biggr]\nonumber
\end{align}
where $p$, $q$ and $r$ vary over all non-negative integers, $f(p) = p$ for $p$ positive and $f(0) = 1$.
\end{theorem}
Further recursion relations, known as the {\em string} and {\em dilaton} equations are satisfied by $\N_{g,n}(\bb)$. See Section~\ref{sec:string-dilaton}.

\begin{figure}[ht] 
\begin{center}
\includegraphics[scale=.9]{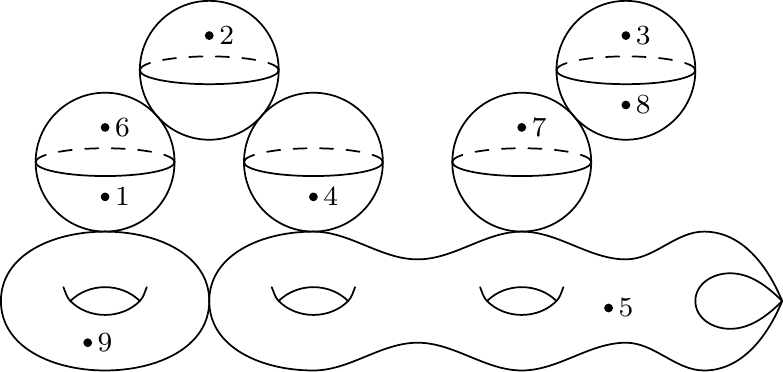} \qquad \qquad \includegraphics[scale=0.8]{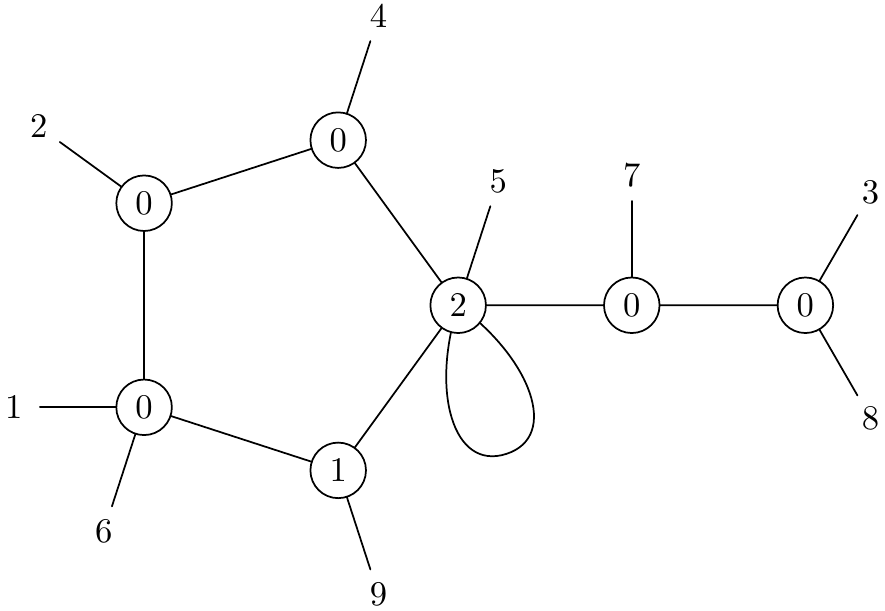}
\caption{Dual graph of a stable curve.}
	\label{fig:dual}
\end{center}
\end{figure}

Recall that the Deligne--Mumford compactification $\M_{g,n}$ possesses a natural stratification indexed by {\em dual graphs}. The dual graph of $\Sigma\in\M_{g,n}$ has vertices corresponding to the irreducible components of $\Sigma$ and assigned genus, edges corresponding to the nodes of $\Sigma$, and a {\em tail}---an edge with an open end (no vertex)---corresponding to each labeled point of $\Sigma$. Figure~\ref{fig:dual} shows an example and Section~\ref{sec:strat} gives precise definitions.

The following theorem expresses $\N_{g,n}(\bb)$ as a sum over dual graphs of type $(g,n)$. Each dual graph contributes the product of its vertex weights divided by the order of its automorphism group. The weight attached to a vertex $v$ is the quasi-polynomial $N_{h(v), n(v)}(\bb_{I(v)}, \mathbf{0})$, where $h(v)$ is the genus of the vertex, $n(v)$ is the valence of the vertex and $I(v)$ denotes the set of labels on the tails adjacent to $v$.

\begin{theorem} \label{th:modsum}
In the following formula the sum is over all dual graphs $G$ of type $(g,n)$ and the product is over the vertices of $G$.
\begin{equation} \label{eq:statesum}
\N_{g,n}(\bb) = \sum_G \frac{1}{|\text{Aut } G|} \prod_{v \in V(G)} N_{h(v), n(v)}(\bb_{I(v)}, \mathbf{0})
\end{equation}
\end{theorem}

\begin{remark} \label{th:drop}
A more natural enumerative problem would be to drop conditions (C3) and (C3') to define $\y_{g,n}(b_1, b_2, \ldots, b_n)\supset\z_{g,n}(b_1, b_2, \ldots, b_n)$ and $\overline{\y}_{g,n}(b_1, b_2, \ldots, b_n)\supset\overline{\z}_{g,n}(b_1, b_2, \ldots, b_n)$ with analogous weighted sums $T_{g,n}(b_1, b_2, \ldots, b_n)$ and $\overline{T}_{g,n}(b_1, b_2, \ldots, b_n)$. In fact $T_{g,n}(b_1, b_2, \ldots, b_n)$ and $\overline{T}_{g,n}(b_1, b_2, \ldots, b_n)$ are determined by and determine $N_{g,n}(b_1, b_2, \ldots, b_n)$ and $\N_{g,n}(b_1, b_2, \ldots, b_n)$. Analogues of Theorems~\ref{th:recursion} and \ref{th:modsum} still hold for $T_{g,n}(b_1, b_2, \ldots, b_n)$ and $\overline{T}_{g,n}(b_1, b_2, \ldots, b_n)$ however their dependence on $b_1, b_2, \ldots, b_n$ is no longer quasi-polynomial and they are more difficult to calculate.
\end{remark}
\begin{remark}
The space $\overline{\z}_{g,n}(b_1, b_2, \ldots, b_n)$ is naturally a suborbifold of the moduli space of stable maps $\M_{g,n}(\PP^1,d)$ for $d = b_1 + b_2 + \cdots + b_n$. Moreover, $\N_{g,n}(\bb)$ (virtually) counts all stable maps satisfying the constraints (C1), (C2) and (C3'). This is not {\em a priori} clear because there are stable maps with domains that are not stable curves. The stable maps with unstable domain have domain with a genus zero irreducible component that maps onto $\bp^1$ and hence has exactly one labeled point (the pre-image of $\infty$) and one node. They contribute a factor of $N_{0,2}(b,0)=0$ by an extension of Theorem~\ref{th:modsum} from stable curves to nodal curves and hence can be ignored. (Note that the constraints (C1) and (C2) do not exclude stable maps since $T_{0,2}(b,0)$ defined in Remark~\ref{th:drop} does not vanish.) There are difficulties in understanding $\N_{g,n}(\bb)$ in terms of intersection theory in $\M_{g,n}(\PP^1,d)$ and Gromov-Witten invariants. One difficulty is that different components of $\overline{\z}_{g,n}(b_1, b_2, \ldots, b_n)$ occur with different multiplicities in $\M_{g,n}(\PP^1,d)$. Another difficulty is relating the virtual count---which takes the Euler characteristic of components---to virtual classes that appear in Gromov-Witten theory.
\end{remark}

The structure of the paper is as follows. Theorems~\ref{th:nbarprop} and \ref{th:recursion} use Theorem~\ref{th:modsum}. Section~\ref{sec:fatgraph} contains preparatory material. The proofs of Theorems~\ref{th:nbarprop} and \ref{th:modsum} are contained in Section~\ref{sec:strat}. The proof of Theorem~\ref{th:recursion} is contained in Section~\ref{sec:recursion}. In Section~\ref{sec:euler} we describe recursions between $\chi(\M_{g,n})$.

\section{Stable fatgraphs} \label{sec:fatgraph}

The main tool we use to enumerate smooth curves equipped with a morphism $f:\Sigma\to\bp^1$ satisfying (C1), (C2) and (C3) its fatgraph, also known as ribbon graph or dessin d'enfant, given by $\Gamma=f^{-1}[0,1]\subset\Sigma$. A fatgraph is an isotopy class of embeddings of a graph into an orientable surface with boundary that defines a homotopy equivalence. In this paper a graph may be disconnected, however it may not contain isolated vertices. The {\em length} of a graph is its number of edges. More formally a fatgraph is defined without reference to a surface.
\begin{definition}
A {\em fatgraph} is a graph $\Gamma$ endowed with a cyclic ordering of half-edges at each vertex. It is uniquely determined by the triple $(X,\tau_0,\tau_1)$ where $X$ is the set of half-edges of $\Gamma$---so each edge of $\Gamma$ appears in $X$ twice---$\tau_1:X\to X$ is the involution that swaps the two half-edges of each edge and $\tau_0:X\to X$ the automorphism that permutes cyclically the half-edges with a common vertex. The underlying graph $\Gamma$ has vertices $X_0=X/\tau_0$, edges $X_1=X/\tau_1$ and boundary components $X_2=X/\tau_2$ for $\tau_2=\tau_0\tau_1$. 
\end{definition}
An {\em automorphism} of a fatgraph $\Gamma$ is a permutation $\phi:X\to X$ that commutes with $\tau_0$ and $\tau_1$. It descends to an automorphism of the underlying graph. If $\Gamma$ is connected, the group generated by $\tau_0$ and $\tau_1$ acts transitively on $X$. Thus an automorphism that fixes a half-edge is necessarily trivial since $\phi(E)=E$ implies $\phi(\tau_0 E)=\tau_0 E$ and $\phi(\tau_1 E)=\tau_1 E$.

A fatgraph structure allows one to uniquely thicken the graph to a surface with boundary. In particular it acquires a type $(g,n)$ for $g$ the genus and $n$ the number of boundary components. The following diagram shows a fatgraph of type $(1,1)$ as well as the surface obtained by thickening the graph. The cyclic ordering of the half-edges with a common vertex is induced by the orientation of the page.

\begin{center}
\includegraphics{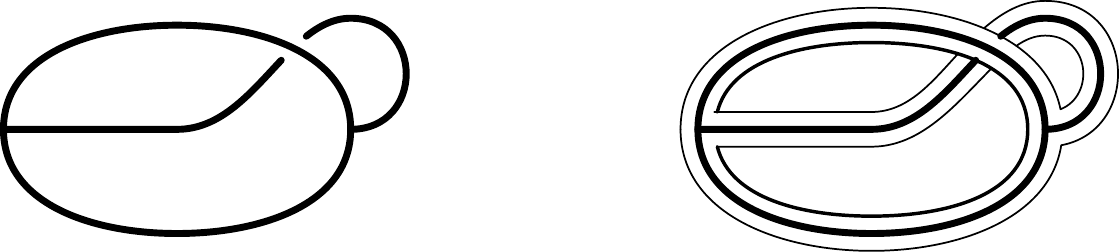}
\end{center}

A {\em labeled} fatgraph is a fatgraph with its boundary components labeled. An {\em automorphism} of a labeled fatgraph $\Gamma$ is a permutation $\phi:X\to X$ that commutes with $\tau_0$ and $\tau_1$ and acts trivially on $X_2$. The automorphism group of a connected labeled fatgraph acts freely on each boundary component since the kernel of the natural restriction map consists of automorphisms that fix a half-edge. In particular it is a subgroup of the rotation group (generated by $\tau_2$) of any boundary component, thus cyclic.

\begin{definition}
For $(b_1, \ldots, b_n)\in\bz_+^n$, define $\fat(b_1, \ldots, b_n)$ to be the set of isomorphism classes of connected, labeled fatgraphs with no valence $1$ vertices, of genus $g$ with $n$ boundary components of lengths $(b_1, \ldots, b_n)$.
\end{definition}

Given a morphism $f:\Sigma\to\bp^1$ satisfying (C1), (C2) and (C3) its fatgraph is given by $\Gamma=f^{-1}[0,1]\subset\Sigma$ with vertices $f^{-1}(0)$ and (centres of) edges $f^{-1}(1)$. Equivalently its set of half-edges $X$ is given by the set of branches of $f^{-1}[0,1]$ with $\tau_0=$ monodromy map around 0 and $\tau_1=$ monodromy map around 1. This defines a map 
\[\z_{g,n}(b_1, b_2, \ldots, b_n)\to\fat(b_1, b_2, \ldots, b_n)\] 
which is an isomorphism. The inverse map is obtained from an explicit construction of a Riemann surface by gluing together $\sum b_i$ copies of $\bc-[0,1]$. The construction also shows that automorphisms of the fatgraph induce automorphisms of the pair $(\Sigma,f)$, so we get \cite{NorCel,NorCou}:
\begin{equation} \label{eq:Nfat}
N_{g,n}(b_1, b_2, \ldots, b_n)=\sum_{\Gamma\in \fat(\bb_S)}\frac{1}{|{\rm Aut\ }\Gamma|}.
\end{equation}
Below we will express $\overline{N}_{g,n}(b_1, b_2, \ldots, b_n)$ as a weighted count of stable fatgraphs.

Kontsevich \cite{KonInt} defined the notion of a stable fatgraph. See also \cite{ZvoStr}.

\begin{definition}
A {\em stable fatgraph} is a fatgraph endowed with the following extra structure.
\begin{itemize}
\item a subset $S\subset X_0$ of distinguished vertices;
\item an equivalence relation $\sim$ on $S$;
\item a genus function $h:S/\hspace{-1mm}\sim\to\bn$ such that $h(S_0)>0$ for any equivalence class $S_0\subset S$ with $|S_0|=1$.
\end{itemize}
\end{definition}
Isomorphisms between stable fatgraphs are isomorphisms of fatgraphs that respect the extra structure---they leave $S$ invariant and preserve $\sim$ and $h$.

Recall that the genus of a connected fatgraph $\Gamma$ is defined by the equation $2-2g-n=V(\Gamma)-E(\Gamma)$ where $V(\Gamma)=|X_0|$, $E(\Gamma)=|X_1|$ and $n=|X_2|$ are the number of vertices, edges and boundary components. More generally, the genus of a connected component $\Gamma'$ of a stable fatgraph $\Gamma$ is defined by removing distinguished vertices so $2-2g(\Gamma')-n(\Gamma')=V(\Gamma'-S)-E(\Gamma')$.
The genus of a stable fatgraph $\Gamma$ requires its {\em dual} graph $G(\Gamma)$. Denote by $\pi_0\Gamma$ the set of connected components of $\Gamma$.
\begin{definition}
Define the {\em dual graph} $G=G(\Gamma)$ of a stable graph $\Gamma$ to have edge set $E(G)=S\cup X_2(\Gamma)$, vertex set $V(G)=(S/\hspace{-1mm}\sim)\cup\pi_0\Gamma$ and incidence relations defined by inclusion. Extend the genus function $h$ on $S/\hspace{-1mm}\sim$ to $h:V(G)\to\bn$ by using the genus of each connected component of $\Gamma$. 
\end{definition}
The genus of a connected (after identification of vertices by $\sim$) stable fatgraph $\Gamma$ is defined to be 
\[ g(\Gamma)=b_1\left(G(\Gamma)\right)+\sum_{v\in V(G(\Gamma))}h(v)\]
where $b_1(G)$ is the first Betti number of $G$.
\begin{definition}
For $(b_1, \ldots, b_n)\in\bz_+^n$, define $\fats(b_1, \ldots, b_n)$ to be the set of isomorphism classes of labeled stable fatgraphs, connected after identification of vertices by $\sim$, of genus $g$ with $n$ boundary components of lengths $(b_1, \ldots, b_n)$, with all vertices of valence 1 contained in $S$.
\end{definition}

One can associate a stable fatgraph to any morphism from a stable curve $f:\Sigma\to\bp^1$ satisfying (C1), (C2) and (C3') as follows. Let $\Gamma'=f^{-1}[0,1]-\{$ nodes, ghost components $\}\subset\Sigma$. Define $\Gamma$ to be the closure of $\Gamma'$ in the normalisation of $\Sigma$, \ie add vertices to non-compact ends of $\Gamma'$. Let $S=\Gamma-\Gamma'$ and define two vertices in $S$ to be equivalent if they coincide in $\Sigma/\hspace{-1.5mm}\sim$ where $\sim$ means ghost components collapsed. The genus $h$ of an equivalence class in $S$ is the genus of the corresponding collapsed component or zero if there is no corresponding collapsed component (so it is purely a node.) This defines a map
\[\overline{\z}_{g,n}(b_1, b_2, \ldots, b_n)\to\fats(b_1, b_2, \ldots, b_n)\] 
which is no longer one-to-one in general since fibres can be infinite. Nevertheless,
\begin{equation} \label{eq:Nfatc}
\overline{N}_{g,n}(b_1, b_2, \ldots, b_n)=\sum_{\Gamma\in\fats(\bb_S)}w(\Gamma)
\end{equation}
for weight $w(\Gamma)$ involving a product of orbifold Euler characteristics of compactified moduli spaces:
\[ w(\Gamma)=\frac{1}{|{\rm Aut\ } \Gamma|} \prod_{v\in S/\sim} \chi\left(\overline{{\mathcal M}}_{h(v),n(v)}\right)\]
where we have defined $n(S_0)=|S_0|$ for any equivalence class $S_0\subset S$ and $\chi(\overline{{\mathcal M}}_{0,2}):=1$ to simplify notation.

As mentioned in Remark~\ref{th:suborb}, $\z_{g,n}(b_1, b_2, \ldots, b_n)$ and $\overline{\z}_{g,n}(b_1, b_2, \ldots, b_n)$ can be identified with suborbifolds of $\modm_{g,n}$, respectively $\M_{g,n}$. Although convenient, it is not essential for the results in this paper so we will simply describe the key ideas and refer the reader to \cite{MPeRib,NorCel} for details. The proof requires one to show that given $\Sigma$, a morphism $f:\Sigma\to\bp^1$ satisfying (C1), (C2) and (C3) is unique and fixed by automorphisms of $\Sigma$. It relies on a theorem due to Strebel \cite{StrQua} which states that for a smooth genus $g$ curve $\Sigma$ with $n$ labeled points and an $n$-tuple $(b_1, b_2, \ldots, b_n) \in \mathbb{R}_+^n$, there exists a unique holomorphic quadratic differential, a {\em Strebel} differential, on $\Sigma-\{p_1, \ldots, p_n\}$ with closed horizontal trajectories and residues at $(p_1, \ldots, p_n)$ determined by $(b_1, \ldots, b_n)$. Furthermore, if $(b_1, \ldots, b_n)$ are positive integers and $\Sigma$ admits a morphism $f:\Sigma\to\bp^1$ satisfying (C1), (C2) and (C3) then the Strebel differential coincides with the pullback $f^*\omega$ for $\omega$ a holomorphic quadratic differential on $\bc-\{0,1\}$. In particular, the uniqueness of $f^*\omega$ implies the uniqueness of $f$. Furthermore, any automorphism of $\Sigma$ fixes $f^*\omega$, by uniqueness of the Strebel differential, and hence fixes $f$. The analogous result for a stable curve $\Sigma\in\overline{\z}_{g,n}(b_1, b_2, \ldots, b_n)$ uses the generalisation of Strebel differentials to stable curves \cite{ZvoStr} which again coincides with the pullback $f^*\omega$.

Connected components of stable fatgraphs consist of fatgraphs with distinguished vertices. It will be convenient to label such vertices when a component is taken in isolation. Such fatgraphs are called {\em pointed} fatgraphs.

\begin{definition} \label{th:point}
A {\em pointed fatgraph} is a labeled fatgraph with some vertices labeled. A {\em pointed stable fatgraph} is a labeled stable fatgraph with some vertices labeled from $X_0-S$.
\end{definition}
Isomorphisms between pointed fatgraphs are isomorphisms of fatgraphs that preserve labeled vertices.
\begin{definition} \label{th:fatpoint}
Define $\fat(b_1, \ldots, b_p, 0, \ldots, 0)$ (respectively $\fats(b_1, \ldots, b_p, 0, \ldots, 0)$), for positive integers $b_1, \ldots, b_p$, to be the set of isomorphism class of pointed (stable) fatgraphs of genus $g$ with $p$ boundary components of lengths $(b_1, \ldots, b_p)$, $n-p$ labeled vertices, all vertices of valence 1 labeled (or contained in $S$), and connected (after identification of vertices by $\sim$.)
\end{definition}
The following proposition is crucial in the proof of Theorem~\ref{th:modsum} which requires us to consider $N_{g,n}(b_1, \ldots, b_n)$ when some $b_i$ vanish.
\begin{proposition} \label{th:pointed}
When some, but not all, of the $b_i$ vanish $N_{g,n}(b_1, \ldots, b_n)$ is a weighted count of pointed fatgraphs. More precisely, for $p>0$ and $b_1, \ldots, b_p$ positive integers
\begin{equation} \label{eq:pointed}
N_{g,n}(b_1, \ldots, b_p, 0, \ldots, 0) = \hspace{-6mm}\sum_{\Gamma \in\fat(b_1, \ldots, b_p, 0, \ldots, 0)} \frac{1}{|\text{Aut } \Gamma|}.
\end{equation}
\end{proposition}
\begin{proof}
Our main tools are the string and dilaton equations for the uncompactified lattice point count proven in \cite{NorStr}. See Section~\ref{sec:strat} for a generalisaton of these equations to $\N_{g,n}$.
\[
N_{g,n+1}(b_1, b_2, \ldots, b_n,2) = \sum_{k=1}^n \sum_{m=1}^{b_k} \left. mN_{g,n}(b_1, b_2, \ldots, b_n) \right|_{b_k=m} - \frac{1}{2} \sum_{k=1}^n b_k N_{g,n}(b_1, b_2, \ldots, b_n)\quad({\rm string\ equation})
\]
\[
N_{g,n+1}(b_1, b_2, \ldots, b_n,2) - N_{g,n+1}(b_1, b_2, \ldots, b_n,0) = (2g-2+n) N_{g,n}(b_1, b_2, \ldots, b_n)\quad({\rm dilaton\ equation}).
\]
The equations apply to the quasi-polynomials and in particular allow some $b_k=0$. In the string equation if $b_k=0$ then the sum over $m=(1, \ldots, b_k)$ does not appear.

We now prove the proposition by induction on $n-p$. The $n-p=0$ case is immediate by definition. Set $\bb_P=(b_1, \ldots, b_p)$. Substitute the string equation into the dilaton equation to obtain the following.
\begin{align*}
& N_{g,n+1}(\bb_P,0, \ldots, 0) \\
=~& N_{g,n+1}(\bb_P,2,0, \ldots, 0) + (2-2g-n) N_{g,n}(\bb_P,0, \ldots, 0) \\
=~& \sum_{k=1}^p \sum_{m=1}^{b_k} \left. mN_{g,n}(\bb_P,0, \ldots, 0) \right|_{b_k = m} - \frac{1}{2} \sum_{k=1}^p b_k N_{g,n}(\bb_P,0, \ldots, 0) + (2-2g-n) N_{g,n}(\bb_P,0, \ldots, 0) \\
=~& \sum_{k=1}^p \sum_{m=1}^{b_k-1} \left. mN_{g,n}(\bb_P,0, \ldots, 0) \right|_{b_k = m} + \left( \frac{1}{2} \sum_{k=1}^p b_k + 2 - 2g - n \right) N_{g,n}(\bb_P,0, \ldots, 0).
\end{align*}
Similar to above, if $b_k-1=0$ then remove the corresponding sum.

The number of vertices of any fatgraph in $\fat(b_1, \ldots, b_n)$ is $V= 2 - 2g - n+\frac{1}{2} \sum_{k=1}^n b_k$ and in particular constant over $\fat(b_1, \ldots, b_n)$. More generally the number of vertices of any pointed fatgraph in $\fat(b_1, \ldots, b_p, 0, \ldots, 0)$ is also constant, given by $V=2 - 2g - p+\frac{1}{2} \sum_{k=1}^p b_k $. 

Then we can rewrite the equation above as
\[
N_{g,n+1}(\bb_P,0, \ldots, 0) = \sum_{k=1}^p \sum_{m=1}^{b_k-1} \left. mN_{g,n}(\bb_P,0, \ldots, 0) \right|_{b_k = m} + (V+p-n) N_{g,n}(\bb_P,0, \ldots, 0).
\]
Put $\fato_{g,n+1}(\bb_P,0, \ldots, 0) =Z_1 \sqcup Z_2$, where $Z_1$ consists of those fatgraphs where the vertex labeled $n+1$ is of valence 1 and $Z_2$ consists of those fatgraphs where the vertex labeled $n+1$ has valence at least two. We will show that the weighted enumeration of $Z_1$ is equal to the first term on the right hand side above while the weighted enumeration of $Z_2$ is equal to the second term on the right hand side above.

\begin{itemize}
\item Note that every pointed fatgraph with a valence one labeled vertex must have trivial automorphism group (since the half-edge incident to the labeled vertex is fixed.) So what we wish to prove is
\[
|Z_1| = \sum_{k=1}^p \sum_{m=1}^{b_k-1} \left. mN_{g,n}(\bb, \mathbf{0}) \right|_{b_k = m}.
\]
We can construct a fatgraph in $Z_1$ by taking a fatgraph $\left. \Gamma'\in\fat(\bb_P,0, \ldots, 0) \right|_{b_k=m}$ and adding a long edge of length $\frac{1}{2}(b_k-m)$ with a vertex on the end labeled $n+1$. Since $\Gamma'$ is connected, $\text{Aut } \Gamma'$ acts freely on the $k$th boundary, so the number of distinct ways to attach the chain is $m/|\text{Aut } \Gamma'|$. Here $m$ is an integer satisfying $0 < m < b_k$ and the construction works for any $k = 1, 2, \ldots, p$ with $b_k-1>0$.

\noindent Therefore, we have
\begin{align*}
\sum_{\Gamma \in Z_1} \frac{1}{|\text{Aut } \Gamma|} = \sum_{\Gamma \in Z_1} 1 &= \sum_{k=1}^p \sum_{m=1}^{b_k-1} \sum_{\Gamma' \in \left. \fat(\bb_P,0, \ldots, 0) \right|_{b_k=m}} \frac{m}{|\text{Aut } \Gamma'|}\\
& = \sum_{k=1}^p \sum_{m=1}^{b_k-1} m \sum_{\Gamma' \in \left. \fat(\bb_P,0, \ldots, 0) \right|_{b_k=m}} \frac{1}{|\text{Aut } \Gamma'|}\\
& = \sum_{k=1}^p \sum_{m=1}^{b_k-1} \left. mN_{g,n}(\bb_P,0, \ldots, 0)\right|_{b_k = m}
\end{align*}

\item The second term counts fatgraphs of type $(g,n)$ with perimeters $\bb_P$ and $n-p$ labeled vertices, and we wish to label one more vertex from the $V+p-n$ unlabeled vertices. Denote by $V_0(\Gamma')$ the unlabeled vertices of $\Gamma'$ (so $|V_0(\Gamma')|=V+p-n$.)

\begin{align*}
(V+p-n) N_{g,n}(\bb_P,0, \ldots, 0) &= \sum_{\substack{\Gamma' \in \fat(\bb_P,0, \ldots, 0) \\ v \in V_0(\Gamma')}} \frac{1}{|\text{Aut } \Gamma'|} = \sum_{\substack{\Gamma' \in \fat(\bb_P,0, \ldots, 0) \\ v \in V_0(\Gamma') / \text{Aut } \Gamma'}} \frac{|\text{Aut }\Gamma'\cdot v|}{|\text{Aut } \Gamma'|}\\
& = \sum_{\substack{\Gamma' \in \fat(\bb_P,0, \ldots, 0) \\ v \in V_0(\Gamma') / \text{Aut } \Gamma'}} \frac{1}{|(\text{Aut } \Gamma')_v|}
\end{align*}
where $\text{Aut }\Gamma'\cdot v$ is the orbit of the vertex $v$ under $\text{Aut }\Gamma'$ and $(\text{Aut } \Gamma')_v\subset\text{Aut } \Gamma'$ is the isotropy subgroup of $v$. Construct $\Gamma\in\fato_{g,n+1}(\bb_P,0, \ldots, 0)$ by labeling the vertex $v\in\Gamma'$. 
The forgetful map induces the exact sequence $1\to\text{Aut } \Gamma\to\text{Aut } \Gamma'$ and since $\text{Aut } \Gamma$ must fix its $n+1-p$ labeled vertices $(\text{Aut } \Gamma')_v$ is its image, \ie 
$(\text{Aut } \Gamma')_v\cong\text{Aut } \Gamma$. Hence
\[
(V+p-n) N_{g,n}(\bb_P,0, \ldots, 0) = \sum_{\Gamma \in Z_2} \frac{1}{|\text{Aut } \Gamma|}.
\]
\end{itemize}
This accounts for all fatgraphs of type $(g,n+1)$ with perimeters prescribed by $(\bb_P,0, \ldots, 0)$ and $n+1-p$ vertices labeled $p+1, p+2, \ldots, n+1$ and the proposition is proven.
\end{proof}

\section{Stratification of $\M_{g,n}$} \label{sec:strat}

The Deligne--Mumford compactification $\M_{g,n}$ possesses a natural stratification by topological type and labeling. To each stable curve $\Sigma$, we associate a combinatorial structure known as a {\em dual graph}. It is a graph with one vertex for each irreducible component of $\Sigma$. The half-edges adjacent to a vertex in the dual graph correspond to distinguished points --- that is, nodes or labeled points --- on the corresponding irreducible component of $\Sigma$. A node is represented by an edge whose endpoints correspond to the components that meet at the node. A labeled point is represented by a half-edge adjacent to a vertex at only one end --- we call these {\em tails.} Each vertex is assigned the geometric genus\footnote{The {\em geometric genus} of an irreducible curve is the genus of its normalisation.} of the corresponding component while each tail is assigned the label of the corresponding labeled point. This discussion motivates the following more precise definition.

\begin{definition}
A {\em dual graph} of type $(g,n)$ is a connected graph $G$ which has $n$ tails and the following extra structure.
\begin{itemize}
\item A bijection which assigns the labels $\{1, 2, \ldots, n\}$ to the tails.
\item A map $h: V(G) \to \{0, 1 ,2, \ldots \}$ which assigns a genus to each vertex of $G$ such that
\[
g = b_1(G) + \sum_{v \in V(G)} h(v).
\]
Each vertex of genus 0 must have valence at least three and each vertex of genus 1 must have valence at least one.
\end{itemize}
\end{definition}

Two dual graphs are isomorphic if and only if there exists a graph isomorphism between them which preserves the genus of each vertex and the label of each tail. As usual, we refer to an isomorphism from a dual graph to itself as an automorphism.

\begin{example} \label{modgraphs}
Up to isomorphism, there are exactly five dual graphs of type $(1,2)$. These are pictured below and their automorphism groups have orders 1, 2, 1, 2, 2, respectively.
\begin{center}
\includegraphics[scale=0.8]{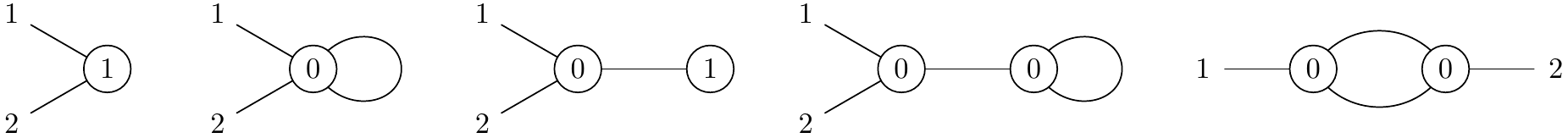}
\end{center}
\end{example}

If $G$ is a dual graph of type $(g,n)$, then the collection of curves ${\mathcal M}_G$ whose associated dual graph is $G$ forms a stratum of $\M_{g,n}$. The stratum ${\mathcal M}_G$ is canonically a product of uncompactified moduli spaces of curves modulo the action of the automorphism group of $G$. Hence, the stratification of $\M_{g,n}$ may be expressed as
\begin{equation} \label{eq:strat}
\M_{g,n} = \bigsqcup_G \prod_{v \in V(G)} {\mathcal M}_{h(v),n(v)} / \textup{Aut } G,
\end{equation}
where the disjoint union is over dual graphs of type $(g,n)$. Here, $n(v)$ denotes the valence of the vertex $v$ while $\textup{Aut } G$ denotes the automorphism group of $G$. Note that there exists a unique open dense stratum formed by the set of smooth curves ${\mathcal M}_{g,n} \subset \M_{g,n}$.

\begin{remark} \label{th:dual}
As one would expect, the dual graphs of a stable curve and a stable fatgraph are related. For a stable curve $\Sigma \in \overline{\z}_{g,n}(b_1, b_2, \ldots, b_n)$, the dual graph of $\Sigma$ is obtained from the natural map $\overline{\z}_{g,n}(b_1, b_2, \ldots, b_n)\to\fats(b_1, b_2, \ldots, b_n)$ by taking the dual graph of the resulting stable fatgraph, removing each valence 2 vertex of genus 0, and identifying its incident edges.
\end{remark}

\begin{proof}[Proof of Theorem~\ref{th:modsum}]
We must express $\N_{g,n}(\bb)$ as a sum over dual graphs of type $(g,n)$. We rewrite (\ref{eq:statesum}) for convenience: 
\[
\N_{g,n}(\bb) = \sum_G \frac{1}{|\textup{Aut } G|} \prod_{v \in V(G)} N_{h(v), n(v)}(\bb_{I(v)}, \mathbf{0}).
\]
Each dual graph contributes the product of its vertex weights divided by the order of its automorphism group. The weight attached to a vertex $v$ is the quasi-polynomial $N_{h(v), n(v)}(\bb_{I(v)}, \mathbf{0})$, where $I(v)$ denotes the set of labels on the tails adjacent to $v$. 

The stratification of $\M_{g,n}$ allows us to decompose $\N_{g,n}(\bb)$ as follows
\[
\N_{g,n}(\bb) = \sum_G N_G(\bb)
\]
for 
\[
N_G(\bb):=\chi\left[\overline{\z}_{g,n}(b_1, b_2, \ldots, b_n)\cap{\mathcal M}_G\right]
\]
Using Remark~\ref{th:dual} we can equivalently interpret $N_G(\bb)$ as a weighted enumeration of stable fatgraphs with perimeters prescribed by $\bb$ whose associated dual graph contracts to $G$. 
\[
N_G(\bb) = \frac{1}{|\textup{Aut } G|} \prod_{v \in V(G)} N_{h(v), n(v)}(\bb_{I(v)}, \mathbf{0}),
\]
where each factor corresponds to choosing a component of the stable integral fatgraph. Furthermore, the correct weight is attached to ghost components since it was proven in \cite{NorCou} that the constant coefficient of $N_{g,n}(b_1, b_2, \ldots, b_n)$ is the orbifold Euler characteristic of ${\mathcal M}_{g,n}$:
\begin{equation} \label{eq:euler}
N_{g,n}(0, 0, \ldots, 0) = \chi({\mathcal M}_{g,n})
\end{equation}
and the sum over all orbifold Euler characteristics of strata of a ghost component gives the orbifold Euler characteristics of the ghost component.

Finally, it is necessary to divide by the number of automorphisms of the dual graph since
\[
|\text{Aut } \Gamma| = |\text{Aut } G(\Gamma)| \prod_{\Gamma_i} |\text{Aut } \Gamma_i|
\]
where the product is over connected components $\Gamma_i$ of $\Gamma$.
\end{proof}

\begin{example}
The formula above allows us to calculate the compactified lattice point count $\N_{g,n}(\bb)$ from the uncompactified lattice point count $N_{g,n}(\bb)$. We can use Example~\ref{modgraphs} to apply this to the case $\N_{1,2}^{(0)}(b_1, b_2)$ as follows, using the expressions for $N_{g,n}(\bb)$ that appear in \cite{NorCou}.
\begin{align*}
\N_{1,2}^{(0)}(b_1, b_2) =~& N_{1,2}^{(0)}(b_1, b_2) + \frac{1}{2}N_{0,4}^{(0)}(b_1, b_2, 0, 0) + N_{0,3}^{(0)}(b_1, b_2, 0) N_{1,1}^{(0)}(0) \\
&+ \frac{1}{2}N_{0,3}^{(0)}(b_1, b_2, 0) N_{0,3}^{(0)}(0, 0, 0) + \frac{1}{2}N_{0,3}^{(0)}(b_1, 0, 0) N_{0,3}^{(0)}(b_2, 0, 0) \\
=~& \frac{1}{384}(b_1^4 + b_2^4 + 2b_1^2b_2^2 -12b_1^2 - 12b_2^2 + 32) + \frac{1}{8}(b_1^2 + b_2^2 - 4) - \frac{1}{12} + \frac{1}{2} + \frac{1}{2} \\
=~& \frac{1}{384} (b_1^4 + b_2^4 + 2b_1^2b_2^2 + 36b_1^2 + 36b_2^2 + 192)
\end{align*}
\end{example}

Theorem~\ref{th:modsum} allows us to deduce properties of the compactified lattice point count $\N_{g,n}(\bb)$ from properties of the uncompactified lattice point count $N_{g,n}(\bb)$.

\begin{proof}[Proof of Theorem~\ref{th:nbarprop}]
This uses the following properties of $N_{g,n}(b_1, b_2, \ldots, b_n)$ proven in \cite{NorCou}.
\begin{itemize}
\item The uncompactified lattice point count $N_{g,n}(\bb)$ is a symmetric quasi-polynomial in $b_1^2, b_2^2, \ldots, b_n^2$ of degree $3g-3+n$ in the sense that it is polynomial on each coset of the sublattice $2\bz^n \subset \bz^n$.
\item If $\alpha_1 + \alpha_2 + \cdots + \alpha_n = 3g-3+n$, then the coefficient of $b_1^{2\alpha_1} b_2^{2\alpha_2} \cdots b_n^{2\alpha_n}$ in $N_{g,n}(b_1, b_2, \ldots, b_n)$ is the following intersection number of psi-classes $\psi_1, \psi_2, \ldots, \psi_n \in H^2(\M_{g,n}; \bq)$.
\begin{equation} \label{eq:int}
\frac{1}{2^{5g-6+2n} \alpha_1! \alpha_2! \cdots \alpha_n!} \int_{\M_{g,n}} \psi_1^{\alpha_1} \psi_2^{\alpha_2} \cdots \psi_n^{\alpha_n}
\end{equation}
\end{itemize}

Theorem~\ref{th:modsum} expresses $\N_{g,n}(\bb)$ as a linear combination of products of uncompactified lattice point polynomials, each of which is quasi-polynomial by the first property above. Therefore, the algebra of quasi-polynomials guarantees that $\N_{g,n}(\bb)$ is a quasi-polynomial in $b_1^2, b_2^2, \ldots, b_n^2$. A quasi-polynomial is symmetric if each polynomial defined on a coset of $2\bz^n$ is invariant under permutations that preserve the coset. In the case of $N_{g,n}(\bb)$, this means that it is symmetric under permutations that preserve the parity of the arguments. The algebra of quasi-polynomials preserves symmetry so it follows that $\N_{g,n}(\bb)$ is symmetric.

By virtue of Theorem~\ref{th:modsum}, we can write $\N_{g,n}(\bb) = N_{g,n}(\bb) + \text{lower order terms}$. This is because the contribution from a stratum is a quasi-polynomial in $b_1^2, b_2^2, \ldots, b_n^2$ with degree equal to the complex dimension of the stratum. Therefore, the degree of $\N_{g,n}(\bb)$ is $3g-3+n$ and Equation~(\ref{eq:int}) implies that the top degree coefficients of $\N_{g,n}(\bb)$ store tautological intersection numbers.

Substitute $b_1 = b_2 = \cdots = b_n = 0$ into Theorem~\ref{th:modsum} and invoke Equation~(\ref{eq:euler}) to deduce that
\[
\N_{g,n}(\mathbf{0}) = \sum_G \frac{1}{|\text{Aut } G|} \prod_{v \in V(G)} N_{h(v), n(v)}(\mathbf{0}) = \sum_G \frac{1}{|\text{Aut } G|} \prod_{v \in V(G)} \chi({\mathcal M}_{h(v), n(v)}) = \chi(\M_{g,n}).
\]
Here, we have used the stratification of $\M_{g,n}$ and the fact that the orbifold Euler characteristic satisfies $\chi(X \times Y) = \chi(X) \chi(Y)$ and $\chi(X \setminus Z) + \chi(Z) = \chi(X)$ for $Z$ a subvariety of $X$. \qedhere
\end{proof}

\begin{remark}
For each dual graph $G$ of type $(g,n)$, define
\[
\N_G(\bb)=\sum_{G' \preceq G}N_G(\bb),
\]
where $G' \preceq G$ if and only if ${\mathcal M}_{G'}$ lies in the closure of ${\mathcal M}_G$. The proof of Theorem~\ref{th:nbarprop} immediately adapts to show that $\N_G(\bb)$ is a quasi-polynomial which satisfies $\N_G(\mathbf{0})=\chi\left(\M_G\right)$. Here, $\M_G$ denotes the closure of the stratum ${\mathcal M}_G \subseteq \M_{g,n}$.
\end{remark}

We are now in a position to generalise Proposition~\ref{th:pointed}.

\begin{corollary} \label{th:pointedc}
For $p>0$ and $b_1, \ldots, b_p$ positive integers
\begin{equation} 
\N_{g,n}(b_1, \ldots, b_p, 0, \ldots, 0) = \sum_{\Gamma \in\fats(b_1, \ldots, b_p, 0, \ldots, 0)} \frac{1}{|\text{Aut } \Gamma|}
\end{equation}
where we recall from Definition~\ref{th:fatpoint} that $\fats(b_1, \ldots, b_p, 0, \ldots, 0)$ consists of pointed stable fatgraphs.
\end{corollary}
\begin{proof}
Put $\bb_P=(b_1, \ldots ,b_p)$. From (\ref{eq:statesum})
\[
\N_{g,n}(b_1, \ldots, b_n) = \sum_G \frac{1}{|\text{Aut } G|} \prod_{v \in V(G)} N_{h(v), n(v)}(\bb_{I(v)}, \mathbf{0})
\]
where $I(v)$ denotes the set of labels on the tails adjacent to $v$. Hence
\[
\N_{g,n}(b_1, \ldots, b_p, 0, \ldots, 0) = \sum_G \frac{1}{|\text{Aut } G|} \prod_{v \in V(G)} N_{h(v), n(v)}(\bb_{I(v)}|_{b_{p+1}= \cdots = b_n = 0}, \mathbf{0})
\]
and Proposition~\ref{th:pointed} tells us that each factor is a weighted count of pointed fatgraphs. Note that the labeled vertices avoid the distinguished vertices in the stable fatgraph so we indeed count elements of $\fats(b_1, \ldots, b_p, 0, \ldots, 0)$.
\end{proof}

\section{Recursion formula} \label{sec:recursion}

In this section we prove the recursion formula of Theorem~\ref{th:recursion} and the string and dilaton equations. We define a {\em long edge} and {\em loop} to be the two graphs consisting of vertices of valence 2 only and a {\em lollipop} to be a loop union a (possible empty) long edge at a valence 3 vertex.

\begin{proof}[Proof of Theorem~\ref{th:recursion}]
We need to prove the recursion (\ref{eq:rec}) which we write again for convenience.
\begin{align*}
\left(\sum_{i=1}^nb_i\right) \N_{g,n}(\bb_S) &= \sum_{i\neq j} \sum_{p+q=b_i+b_j} f(p)q \N_{g,n-1}(p, \bb_{S \setminus \{i,j\}})\\
+ \frac{1}{2}\sum_i&\sum_{p+q+r=b_i} f(p) f(q) r \biggl[\N_{g-1,n+1}(p, q, \bb_{S\setminus \{i\}})
+ \hspace{-5mm}\sum_{\substack{g_1+g_2=g\\I_1 \sqcup I_2 = S\setminus \{i\}}} \hspace{-3mm}\N_{g_1,|I_1|+1}(p, \bb_{I_1}) \N_{g_2,|I_2|+1}(q, \bb_{I_2}) \biggr]\nonumber
\end{align*}
for $S = \{1,2, 3, \ldots, n\}$, $p$, $q$ and $r$ vary over all non-negative integers, $f(p) = p$ if $p$ is positive and $f(0) = 1$.

The strategy of proof is as follows. Construct any $\Gamma\in\fats(\bb_S)$ from smaller fatgraphs by removing from $\Gamma$ a simple subgraph $\gamma$ to get
\[ \Gamma=\Gamma'\cup\gamma.\]
The subgraph $\gamma$ is a long edge or a lollipop which is the simplest subgraph possible so that the remaining fatgraph $\Gamma'$ is legal. There are five cases for removing a long edge or a lollipop from $\Gamma\in\fats(\bb_S)$, shown in Figures~\ref{fig:case1}, \ref{fig:case2}, \ref{fig:case3}, \ref{fig:case4} and \ref{fig:case5}. The broken line signifies $\gamma$, and the remaining stable fatgraph is $\Gamma-\gamma=\Gamma'\in\fato^{\rm stable}_{g',n'}(\bb'_{S'})$ for $(g',n')=(g,n-1)$ or $(g-1,n+1)$ or $\Gamma'=\Gamma_1\sqcup\Gamma_2$ for the pair $\Gamma_i\in\fato^{\rm stable}_{g_i,n_i}(\bb_i)$, $i=1,2$ such that $g_1+g_2=g$ and $n_1+n_2=n+1$. 

In each case, the automorphism groups of $\Gamma'$ and $\Gamma$ act on the construction as follows. The automorphism group of $\Gamma'$ acts on the locations where we attach the ends of $\gamma$. The isotropy subgroup $\ci'\subset{\rm Aut\ }\Gamma'$ is defined to be the subgroup of automorphisms that fix the locations where we attach the ends of $\gamma$. Similarly, the isotropy subgroup $\ci\subset{\rm Aut\ }\Gamma$ is defined to be the subgroup of automorphisms that fix $\gamma$ (and hence the endpoints of $\gamma$.) A simple fact we will use is that $\ci'=\ci$. This is immediate since any automorphism of $\Gamma'$ which fixes the endpoints of $\gamma$ extends to an automorphism of $\Gamma$ which fixes $\gamma$. Conversely any automorphism of $\Gamma$ which fixes $\gamma$ restricts to an automorphism of $\Gamma'$ which fixes the endpoints of $\gamma$. In the simplest case, when $\Gamma$ is connected, $\ci'$ and $\ci$ are both trivial. 

Each fatgraph $\Gamma\in\fats(\bb_S)$ is produced in many ways, one for each long edge and lollipop $\gamma\subset\Gamma$. The number of such $\gamma$ is not constant over all $\Gamma\in\fats(\bb_S)$ however the weighted count over the lengths of each $\gamma$ is constant since each half-edge of $\Gamma$ can be assigned a unique boundary component so
\[|X|=\sum b_i.\]
We exploit this simple fact by taking each construction $q$ times where $\gamma$ has length $q/2$ so that we end up with $(\sum b_i)$ copies of $\Gamma$, if ${\rm Aut\ }\Gamma$ is trivial. More generally, we will explain in each case how to end up with $(\sum b_i)/|{\rm Aut\ }\Gamma|$ copies of $\Gamma$ which is a summand of $(\sum b_i)\cdot\N_{g,n}(\bb_S)$, the left hand side of (\ref{eq:rec}). 

\fbox{\em Case 1} Choose a fatgraph $\Gamma'\in\fato^{\rm stable}_{g,n-1}(p, \bb_{S \setminus \{i,j\}})$ and in Case 1a add a long edge of length $q/2$ inside the boundary of length $p$ so that $p+q=b_i+b_j$ as in the first diagram in Figure~\ref{fig:case1}. 

\begin{figure}[ht] 
	\centerline{\includegraphics[scale=0.3]{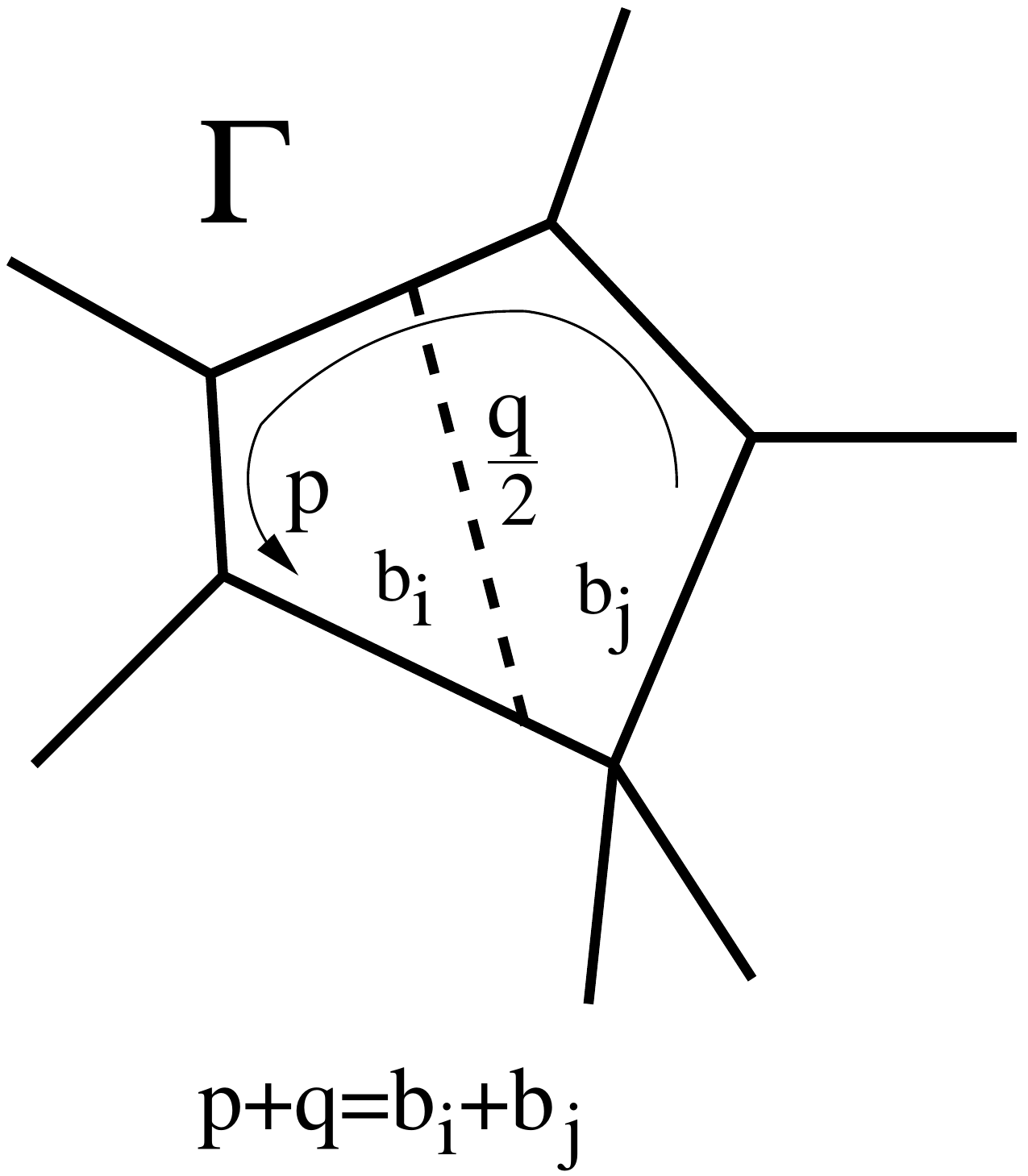} \qquad \qquad \includegraphics[scale=0.3]{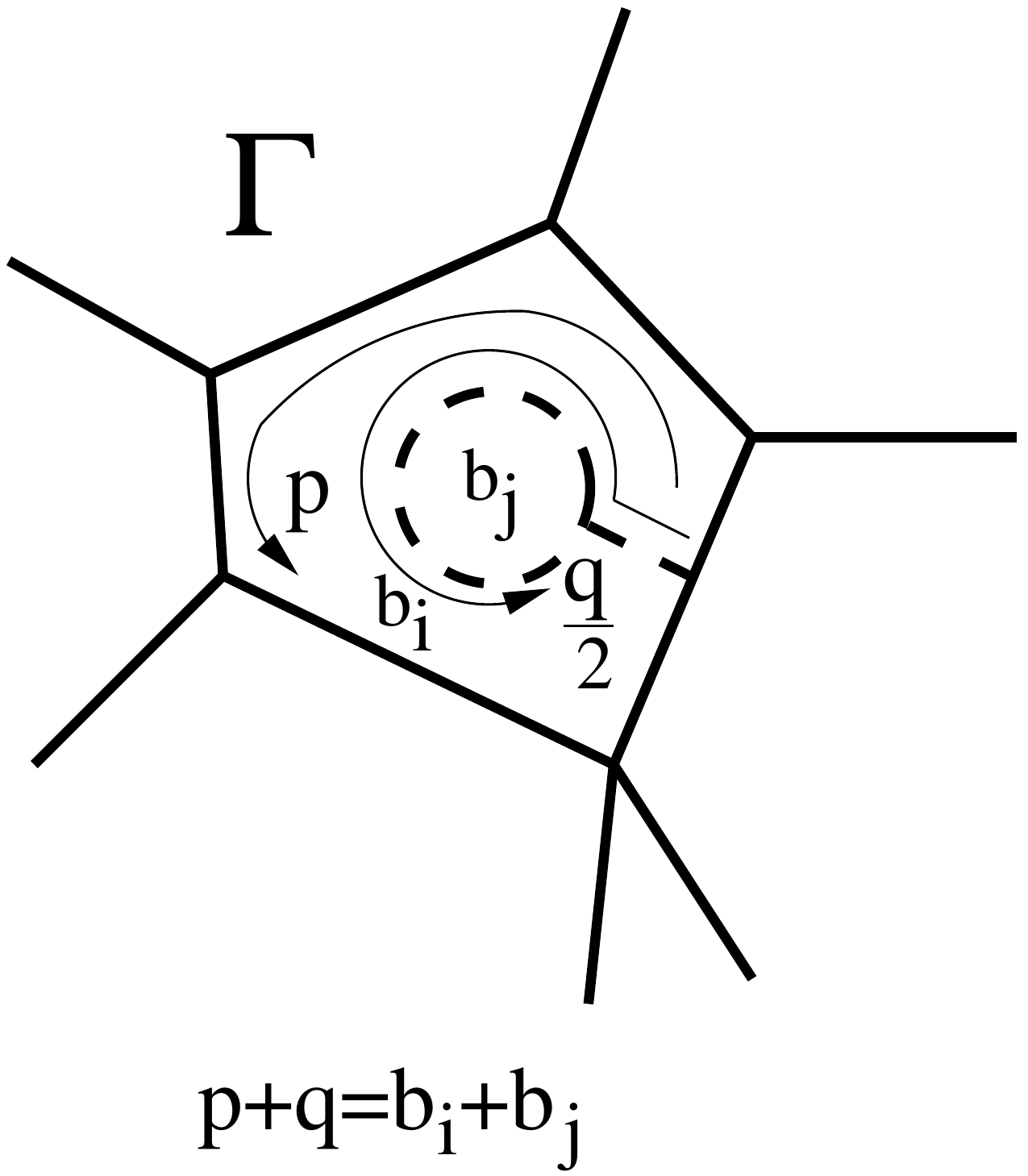}}
	\caption{a. attach edge; b. attach lollipop; to form $\Gamma$.}
	\label{fig:case1}
\end{figure}

In Case 1b attach a lollipop of total length $q/2$ inside the boundary of length $p$ as in the second diagram in Figure~\ref{fig:case1}, again so that $p+q=b_i+b_j$. In both cases for each $\Gamma'$ there are $p$ possible ways to attach the edge, and since the automorphism group of $\Gamma'$ acts on the location where we attach the edge, $q$ copies of this construction produces $pq\cdot|\ci'|/|{\rm Aut\ }\Gamma'|$ stable fatgraphs, where we recall from above that $\ci'\subset{\rm Aut\ }\Gamma'$ is defined to be the subgroup of automorphisms that fix the locations where we attach the ends of $\gamma$. For each $\Gamma$ produced from $\Gamma'$ in this way, this construction produces $q\cdot|\ci|/|{\rm Aut\ }\Gamma|$ copies of $\Gamma$, where again we recall from above that $\ci\subset{\rm Aut\ }\Gamma$ is defined to be the subgroup of automorphisms that fix $\gamma$. Divide by $|\ci'|=|\ci|$ so that $pq/|{\rm Aut\ }\Gamma'|$ stable fatgraphs produce $q/|{\rm Aut\ }\Gamma|$ copies of each $\Gamma$ produced from $\Gamma'$ in this way. Applying this to all $\Gamma'\in\fato^{\rm stable}_{g,n-1}$ this construction contributes 
 \[pq\N_{g,n-1}\left(p, \bb_{S \setminus \{i,j\}}\right)\] 
 to the right hand side of the recursion formula (\ref{eq:rec}) which agrees with a summand.

\fbox{\em Case 2} Choose a pointed fatgraph $\Gamma'\in\fato^{\rm stable}_{g,n-1}(0, \bb_{S \setminus \{i,j\}})$. 
 \begin{figure}[ht] 
	\centerline{\includegraphics[scale=0.3]{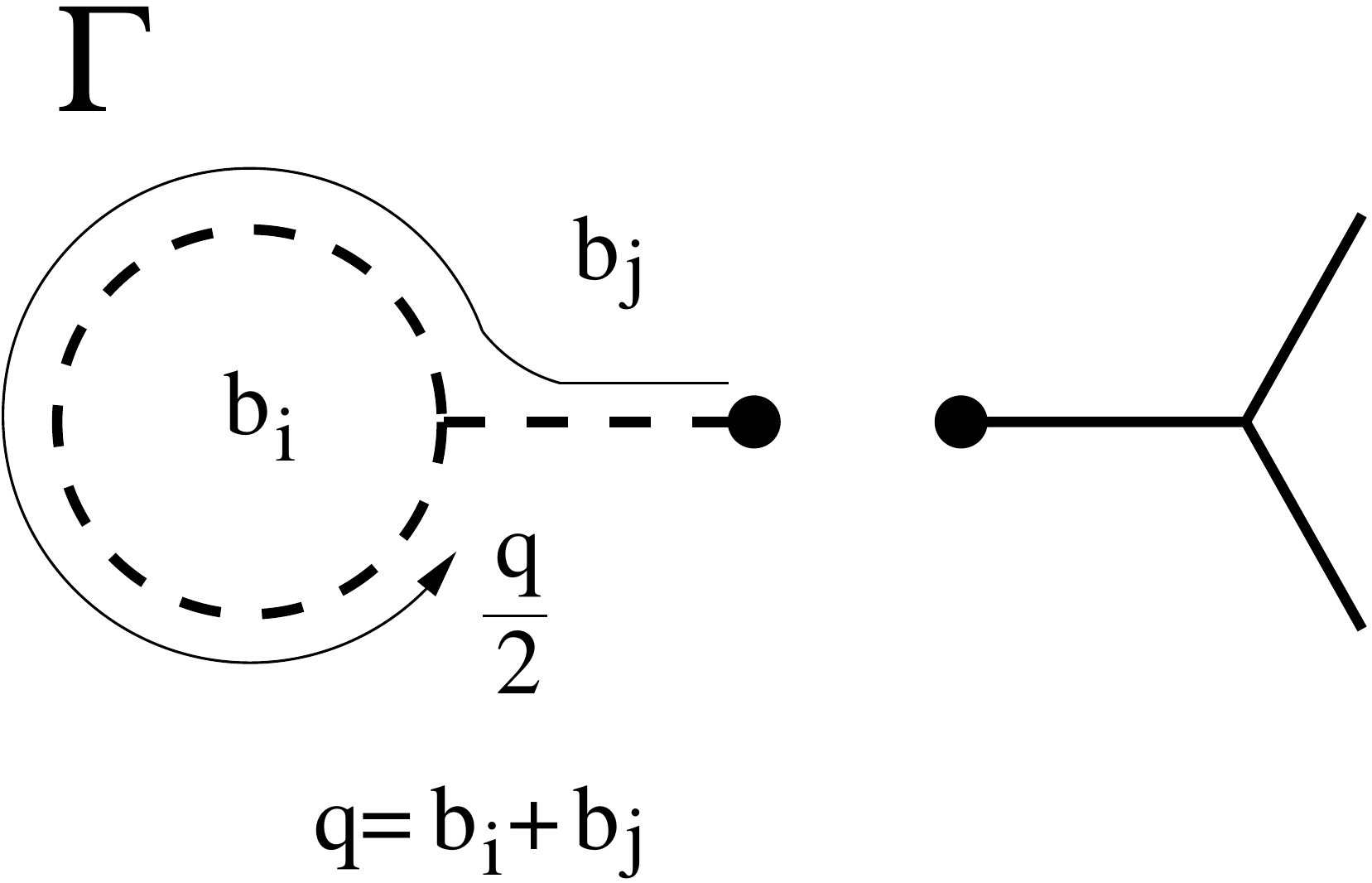}}
	\caption{Identify the vertex of a lollipop with a vertex of $\Gamma'$ to form $\Gamma$.}
	\label{fig:case2}
\end{figure}
Construct $\Gamma$ by identifying the distinguished vertex of $\Gamma'$ with a distinguished vertex of a lollipop. The automorphism group of $\Gamma'$ acts trivially on this construction (since by definition it fixes distinguished vertices) \ie $\ci'={\rm Aut\ }\Gamma'$, so $q$ copies of this construction produces $q$ stable fatgraphs. For each $\Gamma$ produced from $\Gamma'$ in this way, this construction produces $q\cdot|\ci|/|{\rm Aut\ }\Gamma|$ copies of $\Gamma$. Divide by $|\ci'|=|\ci|=|{\rm Aut\ }\Gamma'|$ so that $q/|{\rm Aut\ }\Gamma'|$ stable fatgraphs produce $q/|{\rm Aut\ }\Gamma|$ copies of each $\Gamma$ produced from $\Gamma'$ in this way. Applying this to all $\Gamma'\in\fato^{\rm stable}_{g,n-1}(0, \bb_{S \setminus \{i,j\}})$, and recalling from Corollary~\ref{th:pointedc} that setting a variable to zero counts pointed stable fatgraphs, this construction contributes 
 \[q\N_{g,n-1}\left(0, \bb_{S \setminus \{i,j\}}\right)\] 
 to the right hand side of the recursion formula (\ref{eq:rec}) which agrees with a summand.
This is in some sense a degenerate case of Case 1b, although the pictures show that there is a fundamental difference.

\fbox{\em Case 3} 
Choose a fatgraph $\Gamma'\in\fato^{\rm stable}_{g-1,n+1}(p, q, \bb_{S \setminus \{i\}})$ {\em or} $\Gamma'=\Gamma_1\sqcup\Gamma_2$ for $\Gamma_1\in\fato^{\rm stable}_{g_1,|I_1|+1}(p,\bb_{I_1})$ and $\Gamma_2\in\fato^{\rm stable}_{g_2,|I_2|+1}(q,\bb_{I_2})$ where $g_1+g_2=g$ and $I_1 \sqcup I_2 = S\setminus \{i\}$. Attach a long edge of length $r/2$ connecting these two boundary components as in Figure~\ref{fig:case3} so that $p+q+r=b_i$. 
\begin{figure}[ht] 
	\centerline{\includegraphics[height=4cm]{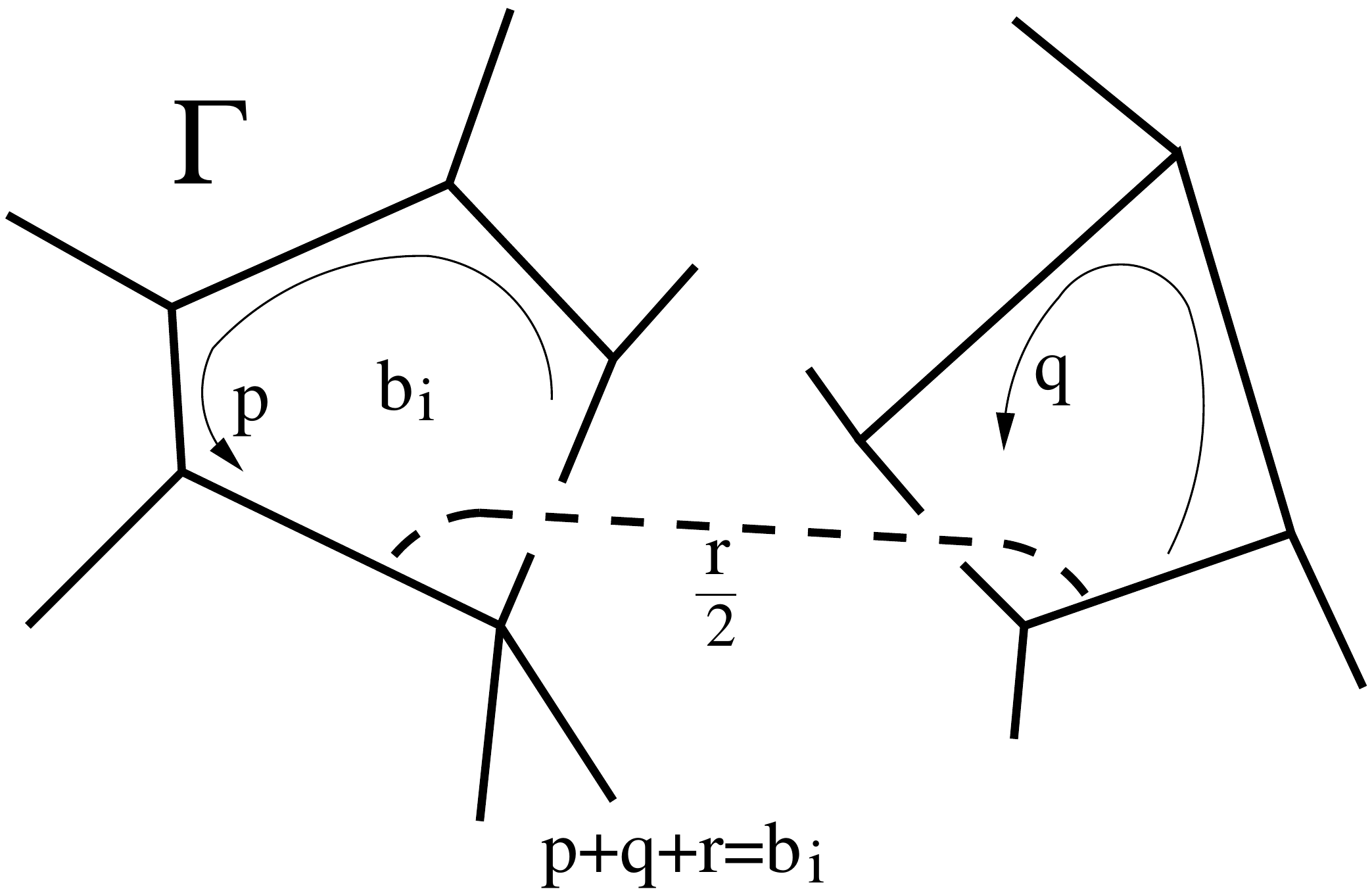}}
	\caption{$\Gamma$ is obtained from a single fatgraph or two disjoint fatgraphs by adding the long edge.}
	\label{fig:case3}
\end{figure}

In the diagram, the two boundary components of lengths $p$ and $q$ are part of a fatgraph that may or may not be connected. There are $pq$ possible ways to attach the edge. An enlarged group of isomorphisms between fatgraphs $\Gamma'$ that does not necessarily preserve the labeling of the two attaching boundary components acts here because we can swap the role of the two attaching boundary components. This either identifies two different fatgraphs $\Gamma'$ or produces new automorphisms of $\Gamma'$. In the first case we count only one of them, or more conveniently we count both of them with a weight of $\frac{1}{2}$. Hence $r$ copies of this construction produces $\frac{1}{2}pqr\cdot|\ci'|/|{\rm Aut\ }\Gamma'|$ stable fatgraphs. In the second case, the action of the automorphism group of $\Gamma'$ on the locations where we attach the edges extends to an action of a larger group ${\rm Aut}^*\Gamma'$ that does not necessarily preserve the labeling of the two attaching boundary components so ${\rm Aut\ }\Gamma'$ is an index 2 subgroup of ${\rm Aut}^*\Gamma'$:
\begin{equation} \label{eq:index2}
1\to{\rm Aut\ }\Gamma'\to{\rm Aut}^*\Gamma'\to\bz_2\to1.
\end{equation}
Hence $r$ copies of this construction produces $pqr\cdot|\ci'|/|{\rm Aut}^*\Gamma'|=\frac{1}{2}pqr\cdot|\ci'|/|{\rm Aut\ }\Gamma'|$ stable fatgraphs so we again count with a weight of $\frac{1}{2}$ as above. For each $\Gamma$ produced from $\Gamma'$ in this way, this construction produces $r\cdot|\ci|/|{\rm Aut\ }\Gamma|$ copies of $\Gamma$. Divide by $|\ci'|=|\ci|$ so that $\frac{1}{2}pqr/|{\rm Aut\ }\Gamma'|$ stable fatgraphs produce $r/|{\rm Aut\ }\Gamma|$ copies of each $\Gamma$ produced from $\Gamma'$ in this way. Applying this to all $\Gamma'\in\fato^{\rm stable}_{g-1,n+1}(p, q, \bb_{S\setminus \{i\}})$ and $\Gamma'=\Gamma_1\sqcup\Gamma_2$ for all $\Gamma_1\in\fato^{\rm stable}_{g_1,j}(p, \bb_{I_1})$ and $\Gamma_2\in\fato^{\rm stable}_{g_2,n+1-j}(q, \bb_{I_2})$
this construction contributes 
\[\frac{1}{2}pqr\biggl[\N_{g-1,n+1}(p, q, \bb_{S\setminus \{i\}})
+ \hspace{-5mm}\sum_{\substack{g_1+g_2=g\\I_1 \sqcup I_2 = S\setminus \{i\}}} \hspace{-3mm}\N_{g_1,|I_1|+1}(p, \bb_{I_1}) \N_{g_2,|I_2|+1}(q, \bb_{I_2})\biggr]\] 
 to the right hand side of the recursion formula (\ref{eq:rec}) which agrees with a summand.

\fbox{\em Case 4} 
Choose a pointed fatgraph $\Gamma'\in\fato^{\rm stable}_{g-1,n+1}(0, q, \bb_{S \setminus \{i\}})$ {\em or} $\Gamma'=\Gamma_1\sqcup\Gamma_2$ for $\Gamma_1\in\fato^{\rm stable}_{g_1,|I_1|+1}(0,\bb_{I_1})$ and 
\begin{figure}[ht] 
 \begin{center}
\includegraphics[scale=0.3]{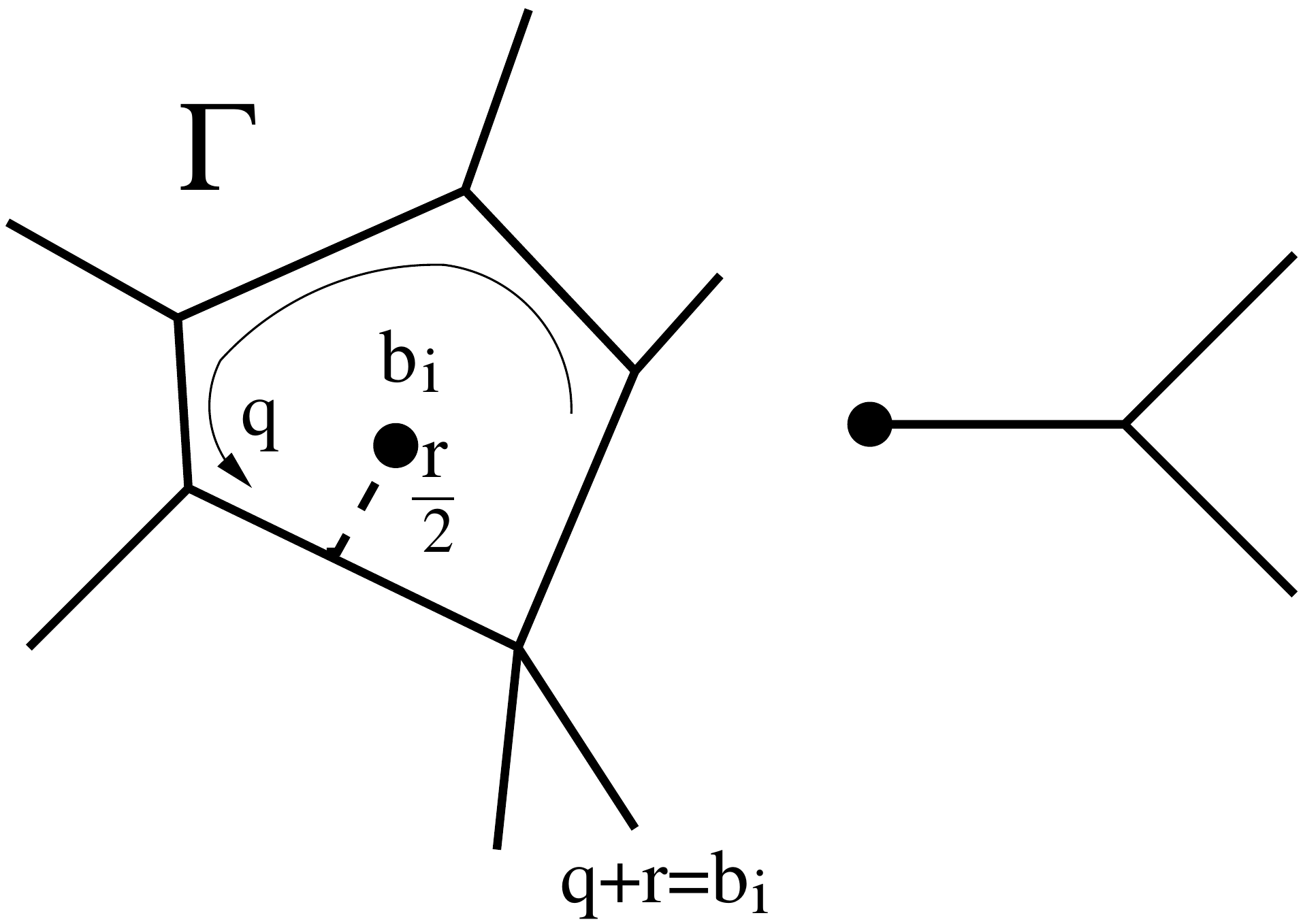}
\end{center}
\caption{$\Gamma$ is obtained from a single fatgraph or two disjoint fatgraphs by adding the broken edge and identifying vertices.}
	\label{fig:case4}
\end{figure}
$\Gamma_2\in\fato^{\rm stable}_{g_2,|I_2|+1}(q,\bb_{I_2})$ where $g_1+g_2=g$ and $I_1 \sqcup I_2 = S\setminus \{i\}$. Attach to a boundary component of $\Gamma'$ or $\Gamma_2$ a long edge of length $r/2$ with a distinguished vertex as in Figure~\ref{fig:case4} so that $q+r=b_i$. 

There are $q$ possible ways to attach the edge, and since the automorphism group of $\Gamma'$ acts on the locations where we attach the edges, $r$ copies of this construction produces $qr|\ci'|/|{\rm Aut\ }\Gamma'|$ stable fatgraphs. For each $\Gamma$ produced from $\Gamma'$ in this way, this construction produces $r\cdot|\ci|/|{\rm Aut\ }\Gamma|$ copies of $\Gamma$. Divide by $|\ci'|=|\ci|$ so that $qr/|{\rm Aut\ }\Gamma'|$ stable fatgraphs produce $r/|{\rm Aut\ }\Gamma|$ copies of each $\Gamma$ produced from $\Gamma'$ in this way. Applying this to all $\Gamma'\in\fato^{\rm stable}_{g-1,n+1}(0, q, \bb_{S\setminus \{i\}})$ and $\Gamma'=\Gamma_1\sqcup\Gamma_2$ for all $\Gamma_1\in\fato^{\rm stable}_{g_1,j}(0, \bb_{I_1})$ and $\Gamma_2\in\fato^{\rm stable}_{g_2,n+1-j}(q, \bb_{I_2})$ this construction contributes 
\[qr\biggl[\N_{g-1,n+1}(0, q, \bb_{S\setminus \{i\}})
+ \hspace{-5mm}\sum_{\substack{g_1+g_2=g\\I_1 \sqcup I_2 = S\setminus \{i\}}} \hspace{-3mm}\N_{g_1,|I_1|+1}(0, \bb_{I_1}) \N_{g_2,|I_2|+1}(q, \bb_{I_2})\biggr]\] 
to the right hand side of the recursion formula (\ref{eq:rec}). It appears with a factor of $\frac{1}{2}$ because (\ref{eq:rec}) includes the isomorphic case of $q=0$ and $p\neq 0$. We have again appealed to Corollary~\ref{th:pointedc} which enables us to count pointed fatgraphs using $\N_{g',n'}(\bb)$ with one of the $b_i=0$.

\fbox{\em Case 5} 
Choose a pointed fatgraph $\Gamma'\in\fato^{\rm stable}_{g-1,n+1}(0, 0, \bb_{S \setminus \{i\}})$ {\em or} $\Gamma'=\Gamma_1\sqcup\Gamma_2$ for $\Gamma_1\in\fato^{\rm stable}_{g_1,|I_1|+1}(0,\bb_{I_1})$ and $\Gamma_2\in\fato^{\rm stable}_{g_2,|I_2|+1}(0,\bb_{I_2})$ where $g_1+g_2=g$ and $I_1 \sqcup I_2 = S\setminus \{i\}$. Identify the two distinguished vertices of a long edge with the two distinguished vertices of $\Gamma'$ as in Figure~\ref{fig:case5} so that $r=b_i$. 
\begin{figure}[ht] 
	\centerline{\includegraphics[height=3cm]{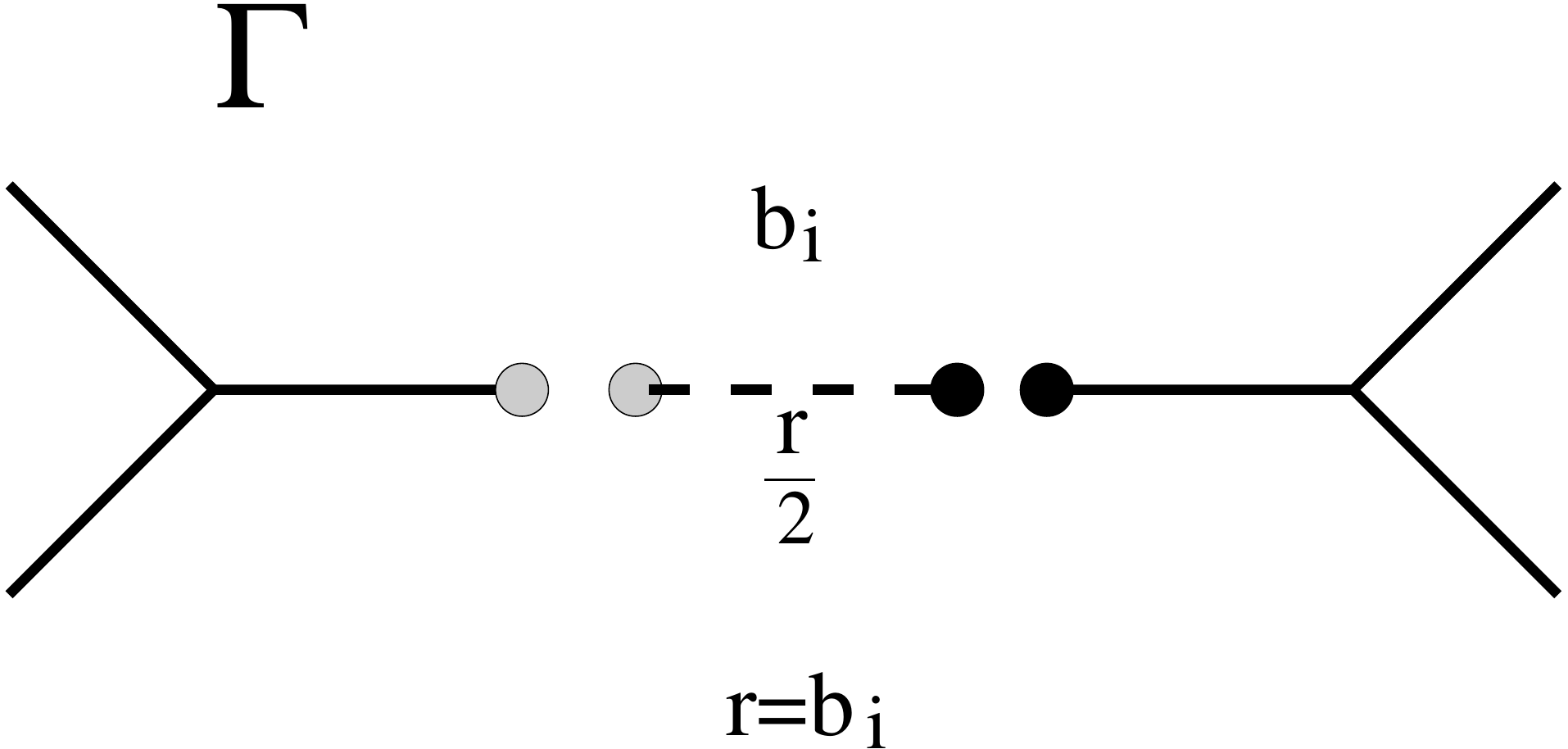}}
	\caption{$\Gamma$ is obtained from a single fatgraph or two disjoint fatgraphs by adding the long edge and identifying vertices.}
	\label{fig:case5}
\end{figure}

In the diagram, the two distinguished vertices are part of a fatgraph that may or may not be connected. The automorphism group of $\Gamma'$ acts trivially on this construction since it fixes distinguished vertices \ie $\ci'={\rm Aut\ }\Gamma'$. As above, an enlarged group of isomorphisms between pointed fatgraphs $\Gamma'$ that does not necessarily preserve the labeling of the two distinguished vertices acts here because we can swap the role of the two distinguished vertices. This either identifies two different fatgraphs $\Gamma'$ or produces new automorphisms of $\Gamma'$. In the first case we count both of them with a weight of $\frac{1}{2}$. Hence $r$ copies of this construction produces $\frac{1}{2}r$ stable fatgraphs. In the second case, the action of the automorphism group of $\Gamma'$ on the locations where we attach the edges extends to an action of the larger group ${\rm Aut}^*\Gamma'$ as in (\ref{eq:index2}). Hence $r$ copies of this construction produces $r\cdot|\ci'|/|{\rm Aut}^*\Gamma'|=\frac{1}{2}r\cdot|\ci'|/|{\rm Aut\ }\Gamma'|=\frac{1}{2}r$ stable fatgraphs which produces a weight of $\frac{1}{2}$ as in the first case. For each $\Gamma$ produced from $\Gamma'$ in this way, this construction produces $r\cdot|\ci|/|{\rm Aut\ }\Gamma|$ copies of $\Gamma$. Divide by $|\ci'|=|\ci|=|{\rm Aut\ }\Gamma'|$ so that $\frac{1}{2}r/|{\rm Aut\ }\Gamma'|$ stable fatgraphs produce $r/|{\rm Aut\ }\Gamma|$ copies of each $\Gamma$ produced from $\Gamma'$ in this way. Applying this to all $\Gamma'\in\fato^{\rm stable}_{g-1,n+1}(0, 0, \bb_{S\setminus \{i\}})$ and $\Gamma'=\Gamma_1\sqcup\Gamma_2$ for all $\Gamma_1\in\fato^{\rm stable}_{g_1,j}(0, \bb_{I_1})$ and $\Gamma_2\in\fato^{\rm stable}_{g_2,n+1-j}(0, \bb_{I_2})$
this construction contributes 
\[\frac{1}{2}r\biggl[\N_{g-1,n+1}(0, 0, \bb_{S\setminus \{i\}})
+ \hspace{-5mm}\sum_{\substack{g_1+g_2=g\\I_1 \sqcup I_2 = S\setminus \{i\}}} \hspace{-3mm}\N_{g_1,|I_1|+1}(0, \bb_{I_1}) \N_{g_2,|I_2|+1}(0, \bb_{I_2})\biggr]\] 
 to the right hand side of the recursion formula (\ref{eq:rec}) which agrees with a summand.

By removing any long edge or lollipop from $\Gamma\in\fat^{\rm stable}(\bb_S)$ we see that it can be produced (many times) using the five constructions above. Each construction produces $\Gamma$ weighted by the factor $2|\gamma|/|{\rm Aut\ }\Gamma|$ where $|\gamma|$ is the length of the long edge or lollipop. The sum over $|\gamma|$ for all long edges or lollipops $\gamma\subset\Gamma$ yields the number of edges of $\Gamma$ so using $|X|=\sum b_i$ this gives a weight of $(\sum b_i)/|{\rm Aut\ }\Gamma|$ to each $\Gamma\in\fat^{\rm stable}(\bb_S)$. The weighted sum over all $\Gamma\in\fat^{\rm stable}(\bb_S)$ is thus $(\sum b_i)\N_{g,n}(\bb_S)$ which gives the left hand side of (\ref{eq:rec}) and completes the proof.
\end{proof}

\begin{remark} \label{th:N11}
A similar argument to the proof of Theorem~\ref{th:recursion} can be used to prove 
\[ b\N_{1,1}^{(0)}(b)=\frac{1}{2}\sum_{p+q+p=b}f(p)q=\frac{b}{2}+\frac{1}{2}\sum_{p+q+p=b}pq=\frac{b}{48}(b^2+20).\]
\end{remark}

\begin{example}
Here we use the recursion formula (\ref{eq:rec}) to calculate $\N_{1,2}(2,2)$. The first sum involves terms with $p+q=2+2=4$ (and $q$ even else the summand vanishes) so $(p,q)=(0,4)$ or $(2,2)$. The second sum involves terms with $p+q+r=2$ (and $r$ even else the summand vanishes) so $(p,q,r)=(0,0,2)$ and there are two terms, one for each boundary component. Thus
\[ 4\N_{1,2}(2,2)=4\N_{1,1}(0)+4\N_{1,1}(2)+\frac{1}{2}\cdot 2\N_{0,3}(0,0,2)+\frac{1}{2}\cdot 2\N_{0,3}(0,0,2)=\frac{17}{3}\]
where we have used $N_{0,3}(0,0,2)=1$ and from Remark~\ref{th:N11}, $\N_{1,1}(0)=5/12$, $\N_{1,1}(2)=1/2$. 
\end{example}

\section{String and dilaton equations} \label{sec:string-dilaton}

It was shown in \cite{NorStr} that the multidifferentials
\[
\omega_{g,n} = \frac{\partial}{\partial z_1} \frac{\partial}{\partial z_2} \cdots \frac{\partial}{\partial z_n} \left( \sum_{b_1, b_2, \ldots, b_n = 1}^\infty N_{g,n}(b_1, b_2, \ldots, b_n) ~ z_1^{b_1} z_2^{b_2} \cdots z_n^{b_n} \right) dz_1~dz_2~\cdots~dz_n
\]
satisfy a topological recursion in the sense of Eynard and Orantin \cite{EOrInv}. One consequence is the fact that there exist string and dilaton equations which provide relations between $\omega_{g,n+1}$ and $\omega_{g,n}$. The corresponding relations between $N_{g,n+1}$ and $N_{g,n}$ are the string and dilaton equations used in the proof of Proposition~\ref{th:pointed}. In the following, we prove that analogous results hold for the compactified lattice point count as well. It would be interesting to know whether the compactified lattice point polynomials can be used to define multidifferentials which also satisfy a topological recursion. 

\begin{theorem}[String equation] \label{string-comp}
Let $f(0) = 1$ and $f(p) = p$ for $p$ positive.
\[
\overline{N}_{g,n+1}(b_1, b_2, \ldots, b_n,1) = \sum_{k=1}^n \sum_{m=0}^{b_k} \left. f(m) \overline{N}_{g,n}(b_1, b_2, \ldots, b_n) \right|_{b_k=m}
\]
\end{theorem}

\begin{proof}
If $b_1 + b_2 + \cdots + b_n$ is even, then both sides of the equation should be interpreted as zero, in which case there is nothing to prove. On the other hand, if $b_1 + b_2 + \cdots + b_n$ is odd, then the inner summation on the right hand side yields a non-zero contribution if and only if $m$ has opposite parity to $b_k$. We may write the string equation as
\begin{equation} \label{string}
\N_{g,n+1}(\bb, 1) = \sum_{k=1}^n \sum_{m=1}^{b_k} \left. m\N_{g,n}(\bb) \right|_{b_k=m} + \sum_{k=1}^n \left. \N_{g,n}(\bb) \right|_{b_k=0}.
\end{equation}

Consider $\Gamma \in \f\hspace{-.3mm}{\rm at}^{\rm stable}_{g,n+1}(\bb, 1)$. The boundary with perimeter 1 belongs to a unique lollipop so suppose that the lollipop is surrounded by boundary $k$. If the long edge of the lollipop has length $a$, then we may write $b_k = m + 2a + 1$, where $m$ is the perimeter of the boundary remaining once the lollipop is removed. After removing the lollipop, the remaining fatgraph $\Gamma'$ is either stable and is an element of $\left.\fats(\bb) \right|_{b_k = m}$ or it is unstable. 

In the first case, $\text{Aut } \Gamma'$ acts on the set $V_k$ of vertices around the boundary labeled $k$ and $\text{Aut } \Gamma = (\text{Aut } \Gamma')_v$ is the isotropy subgroup of automorphisms that fix vertex $v$ where we attach the lollipop. Attaching the lollipop at different vertices in the orbit $(\text{Aut } \Gamma')v$ results in the same fatgraph. Therefore, we obtain the following contribution to $\N_{g,n+1}(\bb, 1)$, where the summation is over $\Gamma' \in \left.\fats(\bb) \right|_{b_k = m}$.
\[
\sum_{\Gamma'} \sum_{v \in (\text{Aut } \Gamma')v} \frac{1}{|(\text{Aut } \Gamma')_v|} = \sum_{\Gamma} \sum_{v \in V_k} \frac{1}{|(\text{Aut } \Gamma')_v| \cdot |(\text{Aut } \Gamma')v|} = \sum_{\Gamma'} \sum_{v \in V_k} \frac{1}{|\text{Aut } \Gamma'|} = \sum_{\Gamma'} \frac{m}{|\text{Aut } \Gamma'|}
\]
Summing over the possible values of $k$ and $m$ yields the first term on the right hand side of (\ref{string}).

In the second case, removing the lollipop leaves an unstable fatgraph precisely when the lollipop belongs to a component of $\Gamma$ of type $(0,3)$. Removal of this component leaves a pointed stable fatgraph $\Gamma'$ of type $(g,n)$ with a distinguished vertex where the extra component is to be attached. Note that $\text{Aut } \Gamma' = \text{Aut } \Gamma$ since the new component has trivial automorphism group and does not introduce any new automorphisms of the corresponding dual graph. Therefore, we obtain the following contribution to $\N_{g,n+1}(\bb, 1)$, where the summation is over $\Gamma' \in \left.\fats(\bb) \right|_{b_k = 0}$.
\[
\sum_{\Gamma} \frac{1}{|\text{Aut } \Gamma|}
\]
Summing over the possible values of $k$ yields the second term on the right hand side of (\ref{string}).
\end{proof}

The proof is purely combinatorial --- the same argument can be used to give a combinatorial proof of the string equation in the uncompactified case.

\begin{theorem}[Dilaton equation] \label{dilaton-comp}
\[
\overline{N}_{g,n+1}(b_1, b_2, \ldots, b_n,2) - \overline{N}_{g,n+1}(b_1, b_2, \ldots, b_n,0) = (2g-2+n) \overline{N}_{g,n}(b_1, b_2, \ldots, b_n).
\]
\end{theorem}

\begin{proof}
The proof relies on the stratification (\ref{eq:strat}) of $\M_{g,n}$ and the dilaton equation for the uncompactified lattice point count. Consider the behaviour of the stratification under the forgetful map $\M_{g,n+1}\to\M_{g,n}$ that forgets $p_{n+1}$. There are two cases.
In the first case on removal of $p_{n+1}$ the underlying curve is still stable, which corresponds to removing a tail of a dual graph. In the second case on removal of $p_{n+1}$ the underlying curve is unstable and the point $p_{n+1}$ lies on a genus zero irreducible component with three distinguished points. There are two ways this can happen---the component has two labeled points and a node; the component has one labeled point, $p_{n+1}$, and two nodes. One can contract the unstable irreducible component, but for our purposes this is unnecessary since each stratum of $\M_{g,n}$ can be obtained from the first case of the forgetful map. The dual graphs from the first case are simply obtained by adding a tail with label $n+1$ to any dual graph of type $(g,n)$.

Using (\ref{eq:statesum}) $\overline{N}_{g,n+1}(b_1, b_2, \ldots, b_n,b_{n+1})$ is the sum of products of $N_{g',n'}(\bb')$ and $b_{n+1}$ appears in exactly one factor of each summand. Hence $\overline{N}_{g,n+1}(b_1, b_2, \ldots, b_n,2) - \overline{N}_{g,n+1}(b_1, b_2, \ldots, b_n,0)$ also factorises with each summand having a factor of the form $N_{g',n'}(\bb_I,2)-N_{g',n'}(\bb_I,0)$. In the first case of the forgetful map discussed above we can use the dilaton equation to get
\[ N_{g',n'}(\bb_I,2)-N_{g',n'}(\bb_I,0)=(2g'-2+n'-1)N_{g',n'-1}(\bb_I)\]
where we note that necessarily $n'>1$. In the second case we get 
\[ N_{0,3}(b_1,b_2,2)-N_{0,3}(b_1,b_2,0)=0.\]
Hence only summands arising from the first case of the forgetful map contribute. 

In terms of pictures, one removes from a dual graph $G$ of type $(g,n+1)$ a tail with label $n+1$ incident to a vertex $v_0$ and replaces it with a dual graph $G'$ of type $(g,n)$ weighted by the factor $2h(v_0) - 2 + n(v_0)-1=2h(v) - 2 + n(v)$. Note that the valence $n(v)=n(v_0)-1$ since an edge is removed. The sum of the weights over all the vertices of a dual graph $G'$ of type $(g,n)$ is
\[\sum_{v \in V(G')} 2h(v) - 2 + n(v)=2g-2+n\] 
which can be understood as relating the arithmetic genus of a stable curve to the Euler characteristic of the curve minus its nodes.
We restate this algebraically: 
\begin{align*}
& \N_{g,n+1}(b_1, b_2, \ldots, b_n,2) - \N_{g,n+1}(b_1, b_2, \ldots, b_n,0) \\
=& \sum_G \frac{1}{| \text{Aut } G|} \prod_{v \in V(G)} \left. N_{h(v), n(v)}(\bb_{I(v)}, \mathbf{0}) \right|_{b_{n+1} = 2} - \sum_{G } \frac{1}{|\text{Aut } G|} \prod_{v \in V(G)} \left. N_{h(v), n(v)}(\bb_{I(v)}, \mathbf{0}) \right|_{b_{n+1} = 0} \\
=& \sum_G \frac{1}{|\text{Aut } G|} \left[ \left. N_{h(v_0),n(v_0)}(\bb_{I(v_0)}, \mathbf{0}) \right|_{b_{n+1} = 2} - \left. N_{h(v_0),n(v_0)}(\bb_{I(v_0)}, \mathbf{0}) \right|_{b_{n+1} = 0} \right] \prod_{v \in V(G) \setminus \{v_0\}} N_{h(v), n(v)}(\bb_{I(v)}, \mathbf{0}) \\
=& \sum_G \frac{1}{|\text{Aut } G|} [2h(v_0) - 3 + n(v_0)] N_{h(v_0), n(v_0)-1}(\bb_{I(v_0) \setminus \{0\}}, \mathbf{0}) \prod_{v \in V(G) \setminus \{v_0\}} N_{h(v), n(v)}(\bb_{I(v)}, \mathbf{0}) \\
=& \sum_{G'} \frac{1}{|\text{Aut } G'|} \left[ \sum_{v \in V(G')} 2h(v) - 2 + n(v) \right] \prod_{v \in V(G')} N_{h(v), n(v)}(\bb_{I(v)}, \mathbf{0}) \\
=& (2g-2+n) \sum_{G'} \frac{1}{|\text{Aut } G'|} \prod_{v \in V(G')} N_{h(v), n(v)}(\bb_{I(v)}, \mathbf{0}) \\
=& (2g-2+n)\N_{g,n}(b_1, b_2, \ldots, b_n).
\end{align*}
The sums begin over all dual graphs $G$ of type $(g,n+1)$, and end over all dual graphs $G'$ of type $(g,n)$ since, as discussed above, those of type $(g,n+1)$ with non-zero contribution correspond to graphs $G'$ of type $(g,n)$ (union a tail.) 
\end{proof}

\section{Euler characteristics} \label{sec:euler}
It was proven by Harer and Zagier \cite{HZaEul}, and independently by Penner \cite{PenPer}, that the orbifold Euler characteristic of the moduli space of curves is 
\[
\chi({\mathcal M}_{g,n}) = (-1)^n \frac{(2g+n-3)!\;B_{2g}}{2g(2g-2)!},
\]
where $B_0, B_1, B_2, \ldots$ denotes the sequence of Bernoulli numbers. They calculate $\chi({\mathcal M}_{g,n})$ from $\chi({\mathcal M}_{g,1})$ via the relation
\begin{equation} \label{eq:eulerel}
\chi(\modm_{g,n+1})=(2-2g-n)\chi(\modm_{g,n}).
\end{equation}
This follows from $\chi(\modm_{g,n})=\chi(\textup{Mod}_{g,n})$ together with the exact sequence of mapping class groups
\begin{equation} \label{eq:exact} 
1 \rightarrow \pi_1(C - \{p_1, p_2, \ldots, p_n\}) \rightarrow \textup{Mod}_{g,n+1} \rightarrow \textup{Mod}_{g,n} \rightarrow 1,
\end{equation}
which implies that $\chi(\textup{Mod}_{g,n+1}) = \chi(C-\{p_1, p_2, \ldots, p_n\}) \times \chi(\text{Mod}_{g,n})$.

Equation~(\ref{eq:eulerel}) is also a consequence of the following properties of $N_{g,n}$ which appear in \cite{NorCou, NorStr}.
\begin{itemize}
\item[(P1)] Orbifold Euler characteristic: $N_{g,n}(0, \ldots, 0)=\chi(\modm_{g,n})$
\item[(P2)] Dilaton equation: $N_{g,n+1}(0, \ldots, 0, 2) - N_{g,n+1}(0, \ldots, 0) = (2g-2+n) N_{g,n}(0, \ldots, 0)$
\item[(P3)] Vanishing: $N_{g,n+1}(0, \ldots, 0, 2)=0$ for $2g-2+n > 0$
\end{itemize}

There is no known closed formula for $\chi(\M_{g,n})$. The aim of this section is to use the following three analogous properties of $\N_{g,n}$ to deduce a recursion relation for $\chi(\M_{g,n})$. For convenience, we define $\chi(\M_{0,1}) = 0$ and $\chi(\M_{0,2}) = 1$.
\begin{itemize}
\item[(P1')] Orbifold Euler characteristic: $\N_{g,n}(0, \ldots, 0)=\chi(\M_{g,n})$
\item[(P2')] Dilaton equation:$\N_{g,n+1}(0, \ldots, 0, 2) - \N_{g,n+1}(0, \ldots, 0) = (2g-2+n) \N_{g,n}(0, \ldots, 0)$
\item[(P3')] $\displaystyle \N_{g,n+1}(0, \ldots, 0, 2)=\frac{1}{2}\chi(\M_{g-1,n+2})+\frac{1}{2}\sum_{h=0}^g\sum_{k=0}^n\binom{n}{k} \chi(\M_{h,k+1}) \; \chi(\M_{g-h,n-k+1})$
\end{itemize}
Properties (P1') and (P2') are contained in Theorems~\ref{th:nbarprop} and \ref{dilaton-comp}. Property (P3') is not a vanishing result, so the recursion relation for $\chi(\M_{g,n})$ is necessarily more complicated than Equation~(\ref{eq:eulerel}).

\begin{proof}[Proof of (P3')]
We begin with the proof of property~(P3) for $N_{g,n}$ because it is needed later. By Proposition~\ref{th:pointed}, $N_{g,n+1}(0, \ldots, 0, 2)$ counts pointed fatgraphs consisting of one edge, one boundary component, and $n$ labeled vertices. But a fatgraph with one edge is either a loop, which we ignore since it has two boundary components, or an edge whose endpoints are necessarily valence one labeled vertices. Thus $N_{g,n+1}(0, \ldots, 0,2)=0$ unless $(g,n)=(0,2)$, in which case we have $N_{0,3}(0, 0, 2) = 1$.

To prove property~(P3'), apply Theorem~\ref{th:modsum}
\[
\N_{g,n+1}(0, \ldots, 0, 2) = \sum_G \frac{1}{|\text{Aut } G|} \prod_{v \in V(G)} N_{h(v), n(v)}(\bb_{I(v)}, \mathbf{0})|_{b_1 = \cdots = b_n = 0,\ b_{n+1}=2}
\]
and note that property~(P3) implies that most terms on the right hand side vanish. The only terms that do not vanish involve $N_{0,3}(0,0,2)=1$ and $N_{h(v), n(v)}(0, \ldots, 0)=\chi(\modm_{h(v), n(v)})$. A useful way to understand this is to consider the stable curves associated to these terms. Recall that a pointed stable fatgraphs enumerated by $\N_{g,n+1}(0, \ldots, 0, 2)$ correspond to a genus $g$ stable curve $\Sigma$ with $n+1$ labeled points, equipped with a morphism $\Sigma \to \bp^1$ that violates C2 --- it may send labeled points to $0$.
\begin{figure}[ht] 
	\centerline{\includegraphics{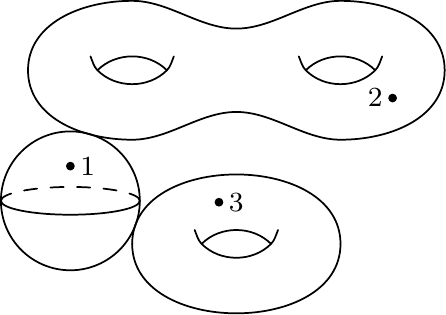}}
	\caption{All components not containing the point labeled 1 map to 0.}
	\label{fig:nodal}
\end{figure}

Denote by $\Sigma_1\subset\Sigma$ the irreducible component containing $p_1$. It necessarily has genus 0. The morphism $\Sigma\to\bp^1$ restricts to a double cover $\Sigma_1\to\bp^1$ ramified at $p_1\mapsto\infty$ and above 1, and $\Sigma-\Sigma_1\mapsto 0$. Each irreducible component of $\Sigma-\Sigma_1$ can vary in its entire moduli space and the weight attached is the Euler characteristic of the moduli space. This suggests how to assemble the different nonvanishing terms---vary a {\em connected} component of $\Sigma-\Sigma_1$ in its {\em compactified} moduli space. Figure~\ref{fig:nodal} shows an example when the complement $\Sigma-\Sigma_1$ is disconnected. The contribution of connected components of arithmetic genus $j$, respectively $g-j$, and $k+1$, respectively $n-k+1$, labeled points to $\N_{g,n+1}(2, 0, \ldots, 0)$ is the weight $\chi(\M_{j,k+1})\chi(\M_{g-j,n-k+1})$. There are $\binom{n}{k}$ ways to partition $n$ labeled points of $\Sigma$ into $k$ and $n-k$ sets. (The nodes account for the extra labeled points.) The factor of $1/2$ appears in front of each summand in (P3') because we either count a decomposition twice or when $n=0$ there exists extra isomorphisms and automorphisms swapping the connected components. If $\Sigma-\Sigma_1$ is connected it has arithmetic genus $g-1$ and $n+2$ labeled points and its contribution is the weight $\frac{1}{2}\chi(\M_{g-1,n+2})$. The factor of $1/2$ appears because the two nodes are in fact not labeled so there exists extra isomorphisms. If $\Sigma-C_1$ consists of a labeled point disjoint union a connected component of arithmetic genus $g$ and $n$ labeled points then its contribution is the weight $n\chi(\M_{g,n})$. We conveniently encode this in (P3') using $\chi(\M_{0,2}):=1$ and including each factor twice, each weighted with a factor of $1/2$. We have accounted for all terms on the right side of (P3') and the equation is proven.
\end{proof}

An immediate consequence of (P1'), (P2') and (P3') is the following analogue of Equation~(\ref{eq:eulerel}).

\begin{proposition} \label{th:euler}
The orbifold Euler characteristics of the compactified moduli spaces of curves satisfy the following recursion relation.
\[
\chi(\M_{g,n+1})=(2-2g-n)\chi(\M_{g,n})+\frac{1}{2}\chi(\M_{g-1,n+2})+\frac{1}{2}\sum_{h=0}^g\sum_{k=0}^n\binom{n}{k}\chi(\M_{h,k+1}) \; \chi(\M_{g-h,n-k+1}).
\]
\end{proposition}

Define the generating function
\[
G(x,q) = \sum_{g=0}^\infty \sum_{n=1}^\infty\frac{\chi(\M_{g,n})}{(n-1)!} x^{n-1} q^g,
\] 
where we take $\chi(\M_{0,1}) = 0$ and $\chi(\M_{0,2}) =  1$. Then Proposition~\ref{th:euler} is equivalent to the PDE
\[
G_x = G + 1-xG_x + GG_x + \frac{q}{2} G_{xx} - 2q G_q.
\]
If we define $F_0(x) = G(x,0) = \sum \frac{\chi(\M_{0,n})}{(n-1)!} x^{n-1}$, we obtain $F_0(0)=0$ and $F_0'=\frac{F_0+1}{1+x-F_0}$. The solution of this differential equation is the inverse of the function
\[
x = 2F_0 - (1+F_0) \ln(1+F_0)
\]
whose expansion is $\displaystyle F_0(x) = x + \frac{1}{2}x^2 + \frac{1}{3}x^3 + \frac{7}{24}x^4 + \frac{17}{60}x^5 + \cdots$. Thus, we recover the genus zero results obtained in \cite{GLSSeq,ManGen}.

The PDE can be studied genus by genus, which yields a hierarchy of linear first order ODEs, each of which contains lower genus solutions. For example, to study the genus one case, we define $F_1(x) = G_q(x,0) = \sum \frac{\chi(\M_{1,n})}{(n-1)!} x^{n-1}$ and we obtain $F_1(0)=\frac{5}{12}$ and
\[
F_1'=\frac{F_1(F_0''-1)+1+F_0''/2}{1+x-F_0}.
\]

\appendix

\section{Table of compactified lattice point polynomials}

\begin{center}
\begin{tabular*}{\textwidth}{@{}cccl@{}} \toprule
$g$ & $n$ & $k$ & $\overline {N}_{g,n}^{(k)}(b_1, b_2, \ldots, b_n)$ \\ \midrule
0 & 3 & 0 & 1 \\ \midrule
0 & 3 & 2 & 1 \\ \midrule
1 & 1 & 0 & $\frac{1}{48}(b_1^2 + 20)$ \\ \midrule
0 & 4 & 0 & $\frac{1}{4}(b_1^2 + b_2^2 + b_3^2 + b_4^2 + 8)$ \\ \midrule
0 & 4 & 2 & $\frac{1}{4}(b_1^2 + b_2^2 + b_3^2 + b_4^2 + 2)$ \\ \midrule
0 & 4 & 4 & $\frac{1}{4}(b_1^2 + b_2^2 + b_3^2 + b_4^2 + 8)$ \\ \midrule
1 & 2 & 0 & $\frac{1}{384} (b_1^4 + b_2^4 + 2b_1^2b_2^2 + 36b_1^2 + 36b_2^2 + 192)$ \\ \midrule
1 & 2 & 2 & $\frac{1}{384} (b_1^4 + b_2^4 + 2b_1^2b_2^2 + 36b_1^2 + 36b_2^2 + 84)$ \\ \midrule
0 & 5 & 0 & $\frac{1}{32} \sum b_i^4 + \frac{1}{8} \sum b_i^2b_j^2 + \frac{7}{8} \sum b_i^2 + 7$ \\ \midrule
0 & 5 & 2 & $\frac{1}{32} \sum b_i^4 + \frac{1}{8} \sum b_i^2b_j^2 + \frac{5}{16} b_1^2 + \frac{5}{16} b_2^2 + \frac{1}{8} b_3^2 + \frac{1}{8} b_4^2 + \frac{1}{8} b_5^2 + \frac{19}{16}$ \\ \midrule
0 & 5 & 4 & $\frac{1}{32} \sum b_i^4 + \frac{1}{8} \sum b_i^2b_j^2 + \frac{5}{16} b_1^2 + \frac{5}{16} b_2^2 + \frac{5}{16} b_3^2 + \frac{5}{16} b_4^2 + \frac{7}{8} b_5^2 + \frac{7}{8}$ \\ \midrule
1 & 3 & 0 & $\frac{1}{4608} \sum b_i^6 + \frac{1}{768} \sum b_i^4b_j^2 + \frac{1}{384} b_1^2b_2^2b_3^2 + \frac{13}{1152} \sum b_i^4 + \frac{1}{24} \sum b_i^2b_j^2 + \frac{29}{144} \sum b_i^2 + \frac{17}{12}$ \\ \midrule
1 & 3 & 2 & $\frac{1}{4608} \sum b_i^6 + \frac{1}{768} \sum b_i^4b_j^2 + \frac{1}{384} b_1^2b_2^2b_3^2 + \frac{43}{4608} (b_1^4 \!+\! b_2^4)\! + \frac{13}{1152} b_3^4 + \frac{1}{24} \sum b_i^2b_j^2 + \frac{277}{4608} (b_1^2 \!+\! b_2^2)\! + \frac{35}{576} b_3^2 + \frac{81}{256}$ \\ \midrule
2 & 1 & 0 & $\frac{1}{1769472}b_1^8 + \frac{3}{40960}b_1^6 + \frac{133}{61440}b_1^4 + \frac{1087}{34560}b_1^2 + \frac{247}{1440}$ \\ \midrule 
0 & 6 & 0 & $\frac{1}{384} \sum b_i^6 + \frac{3}{128} \sum b_i^4b_j^2 + \frac{3}{32} \sum b_i^2b_j^2b_k^2 + \frac{1}{6} \sum b_i^4 + \frac{9}{16} \sum b_i^2b_j^2 + \frac{109}{24} \sum b_i^2 + 34$ \\ \bottomrule
\end{tabular*}
\end{center}


\begin{thebibliography}{99}

\bibitem{BHaEul} Bini, Gilberto and Harer, John.
\emph{Euler characteristics of moduli spaces of curves.} 
J. Eur. Math. Soc. {\bf 13} (2011), 487--512.

\bibitem{EOrInv} Eynard, Bertrand and Orantin, Nicolas.
\emph{Invariants of algebraic curves and topological expansion.}
Commun. Number Theory Phys. {\bf 1} (2007), 347--452.

\bibitem{GLSSeq} Goulden, Ian; Litsyn, Simon and Shevelev, Vladimir.
\emph{On a sequence arising in algebraic geometry.}
J. Integer Seq. {\bf 8} (2005), Article 05.4.7.

\bibitem{HZaEul} Harer, John and Zagier, Don.
\emph{The Euler characteristic of the moduli space of curves.}
Invent. Math. {\bf 85} (1986), 457--485.

\bibitem{KonInt} Kontsevich, Maxim.
\emph{Intersection theory on the moduli space of curves and the matrix Airy function.} 
Comm. Math. Phys. {\bf 147} (1992), 1--23.

\bibitem{ManGen} Manin, Yuri.
\emph{Generating functions in algebraic geometry and sums over trees.} 
The moduli space of curves (Texel Island, 1994), 401--417, Progr. Math., {\bf 129}, Birkh\"auser Boston, Boston, MA, 1995.

\bibitem{MPeRib} Mulase, Motohico and Penkava, Michael.
\emph{Ribbon graphs, quadratic differentials on Riemann surfaces, and algebraic curves defined over $\overline{\mathbb Q}$.}
Mikio Sato: a great Japanese mathematician of the twentieth century.
Asian J. Math. {\bf 2} (1998), 875--919.

\bibitem{NorCel}
Norbury, Paul.
\emph{Cell decompositions of moduli space, lattice points and Hurwitz problems.}
\href{http://arxiv.org/abs/1006.1153}{\texttt{arXiv:1006.1153v1 [math.GT]}}.

\bibitem{NorCou}
Norbury, Paul.
\emph{Counting lattice points in the moduli space of curves.}
Math. Res. Lett. {\bf 17} (2010), 467--481.

\bibitem{NorStr}
Norbury, Paul.
\emph{String and dilaton equations for counting lattice points in the moduli space of curves.}
To appear in Trans. Amer. Math. Soc.

\bibitem{PenPer} Penner, Robert.
\emph{Perturbative series and the moduli space of Riemann surfaces.}
J. Differential Geom. {\bf 27} (1988), 35--53.

\bibitem{StrQua} Strebel, Kurt.
\emph{Quadratic differentials.} 
Springer--Verlag, Berlin, 1984.

\bibitem{WitTwo} Witten, Edward.
\emph{Two-dimensional gravity and intersection theory on moduli space.}
Surveys in differential geometry (Cambridge, MA, 1990), 243--310, Lehigh Univ., Bethlehem, PA, 1991.

\bibitem{ZvoStr} Zvonkine, Dimitri.
\emph{Strebel differentials on stable curves and Kontsevich's proof of Witten's conjecture.} 
\href{http://arxiv.org/abs/math/0209071}{\texttt{arXiv:math/0209071v2 [math.AG]}}.

\end{thebibliography}
\end{document}